\documentclass[10pt,reqno]{amsart}
\usepackage{geometry}
\usepackage{amssymb}
\usepackage{graphicx,subfigure,url}
\usepackage{setspace}

\renewcommand{\Re}{\ensuremath{\mathbb{R}}}

\newcommand{\deriv}[2]{\ensuremath{\frac{\partial #1}{\partial #2}}}
\newcommand{\refsec}[1]{Section \ref{sec:#1}}
\newcommand{\refeqn}[1]{(\ref{eqn:#1})}
\newcommand{\reffig}[1]{Figure \ref{fig:#1}}
\newcommand{\norm}[1]{\ensuremath{\left\| #1 \right\|}}

\newcommand{\bracket}[1]{\ensuremath{\left[ #1 \right]}}
\newcommand{\braces}[1]{\ensuremath{\left\{ #1 \right\}}}
\newcommand{\parenth}[1]{\ensuremath{\left( #1 \right)}}
\newcommand{\tr}[1]{\mbox{tr}\ensuremath{\negthickspace\bracket{#1}}}
\newcommand{\so}{\ensuremath{\mathfrak{so}(3)}}
\newcommand{\SO}{\ensuremath{\mathrm{SO}(3)}}
\newcommand{\SE}{\ensuremath{\mathrm{SE}(3)}}
\renewcommand{\Re}{\ensuremath{\mathbb{R}}}

\renewcommand{\S}{\ensuremath{\mathsf{S}}}
\newcommand{\Le}{\ensuremath{\mathsf{L}}}
\newcommand{\T}{\ensuremath{\mathsf{T}}}
\newcommand{\D}{\ensuremath{\mathbf{D}}}

\title{Computational Geometric Optimal Control of Rigid Bodies}
\author{Taeyoung Lee}
\author{Melvin Leok}
\author{N. Harris McClamroch}

\thanks
{Taeyoung Lee, Department of Aerospace Engineering, University of Michigan. \url{tylee@umich.edu}\endgraf
Melvin Leok, Assistant Professor, Department of Mathematics, Purdue University. \url{mleok@math.purdue.edu}\endgraf
N. Harris McClamroch, Professor, Department of Aerospace Engineering, University of Michigan. \url{nhm@umich.edu}
}

\begin{document}

\maketitle

\begin{abstract}
This paper formulates optimal control problems for rigid bodies in a geometric manner and it presents computational procedures based on this geometric formulation for numerically solving these optimal control problems.   The dynamics of each rigid body is viewed as evolving on a configuration manifold that is a Lie group. Discrete-time dynamics of each rigid body are developed that evolve on the configuration manifold according to a discrete version of Hamilton's principle so that the computations preserve geometric features of the dynamics and guarantee evolution on the configuration manifold; these discrete-time dynamics are referred to as Lie group variational integrators.   Rigid body optimal control problems are formulated as discrete-time optimization problems for discrete Lagrangian/Hamiltonian dynamics, to which standard numerical optimization algorithms can be applied.   This general approach is illustrated by presenting results for several different optimal control problems for a single rigid body and for multiple interacting rigid bodies. The computational advantages of the approach, that arise from correctly modeling the geometry, are discussed.
\end{abstract}

\section{Introduction}

This paper utilizes methods from geometric mechanics and optimal control to develop new computational procedures for geometric optimal control of rigid bodies. The emphasis is on formulating a discrete-time optimal control problem that inherits important conservation properties of rigid body dynamics; this is achieved by combining variational integrators~\cite{MaWe2001} and Lie group methods~\cite{IsMuNoZa2000} to evolve the mechanical configuration. This approach leads to Lie group variational integrators that define the discrete-time rigid body dynamics which the optimal control computations are based upon~\cite{CMA07,CMDA07}.

Most of the prior work related to optimal control of a rigid body is based on local coordinates on $\SO$ or quaternions~\cite{Bil.JoGCaD1993,Bye.JoGCaD1993,Scr.JoGCaD1994,Sey.JoGCaD1993}. Minimal representations of the attitude of a rigid body, such as Euler angles, exhibit coordinate singularities, and require manipulating complicated trigonometric expressions. Nonminimal representations such as quaternions have no coordinate singularities, but they also introduce certain complications. In particular, the group of unit quaternions $\mathsf{SU}(2)\simeq\S^3$ double covers $\SO$, so there is an ambiguity in representing an attitude of a rigid body. Furthermore, the Hamiltonian structure of rigid body attitude dynamics is unnecessarily complicated when it is expressed in terms of quaternions~\cite{LeRe2004}.

By considering rigid body translation and rotation as evolution on a Lie group, optimal control problems defined on Lie groups were introduced by Roger Brockett~\cite{BrockSIJC72,BrockSIAM73} and by John Baillieul~\cite{BailJOTA78}. They emphasized the use of Lie group structures to characterize controllability and existence of optimal controls; they also obtained analytical results for the solution of certain types of optimal control problems. An optimal control problem for a generalized rigid body on $\mathrm{SO}(n)$ was considered in~\cite{Blo.SaCL1996}, and a general theory of optimal control problems on a Lie group was developed in~\cite{Jur.1998a,Jur.1998,Jur.BK1997} together with reachability and controllability conditions. Although these papers viewed rigid body translation and rotation as motion on a Lie group, their results are limited to optimal control problems that can be formulated solely in terms of kinematics. In particular, they do not include dynamics in their analysis, and assume that the controls enter directly at the level of the Lie algebra.

The approach of computational geometric optimal control is focused on developing numerical algorithms, for optimal control problems, that preserve the geometric properties of the dynamics and the optimal control problem~\cite{Lee.2008}. The essential idea is to apply geometric optimal control theory to discrete-time mechanical systems obtained using geometric numerical integrators. A discrete-time version of the generalized rigid body equations and their formulation as an optimal control problem are presented in~\cite{Blo.ICoDaC1998,Blo.N2002}, and discrete-time optimal control problems for the dynamics of a rigid body are considered in~\cite{Blo.JoDaCS2007,Lee.JoDaCS2007,Lee.PotICoDaC2006}. A direct optimal control approach is applied to discrete-time mechanical systems in~\cite{Jun.IC2005}, and it is referred to as \textit{Discrete Mechanics and Optimal Control}.

This paper presents the approach of computational geometric optimal control for the dynamics of rigid bodies on a Lie group. We take the same geometric perspective as in the work of Roger Brockett~\cite{BrockSIJC72,BrockSIAM73}, viewing evolution on a Lie group as fundamental.  However, the emphasis in the present paper is on geometric formulations of both the kinematics and dynamics in the optimal control formulation and the role of geometric methods in optimal control computations.

The development in the paper makes clear that there are important advantages in formulating the optimal control problem as a discrete-time optimal control problem using Lie group variational integrators and then applying standard computational methods to solve the resulting discrete-time optimization problem. This is in contrast with approaches that construct continuous-time necessary conditions and then make use of computational methods to solve these necessary conditions.  The paper demonstrates that for for the optimal control of rigid bodies, the proposed approach exhibits important advantages.

The main contributions of this paper can be summarized as follows: (i) the analytical and computational results presented in this paper are coordinate free; they avoid the singularities, ambiguity, and complications associated with local coordinates, and they provide a global insight into rigid body dynamics, (ii) a geometric optimal control problem is formulated for nontrivial rigid body dynamics that evolve on a Lie group, and (iii) a computational geometric optimal control approach is developed based on a geometric numerical integrator.

\refsec{MFOCRB} provides a summary of Lie group variational integrators for rigid bodies that evolve on a Lie group.   The resulting discrete-time rigid body dynamics are used as a basis for formulating a discrete-time optimal control problem.   In \refsec{optsrbdy} and \ref{sec:optmrbdy}, four different examples of rigid body optimal control problems are studied in some detail.   First, optimal orbit and attitude maneuvers for a rigid dumbbell spacecraft in orbit about a large central body are studied.   Then, optimal attitude maneuvers for a 3D pendulum acting under uniform gravity are studied; the control input conserves the component of the vertical component of the angular momentum thereby requiring a careful computational treatment that avoids numerical ill-conditioning.   The third example is a 3D pendulum attached to a cart that can move in a horizontal plane; optimal reconfiguration maneuvers are studied for this cart and pendulum system.   The fourth example involves optimal attitude maneuvers of two rigid bodies connected by a universal joint;  the control input conserves angular momentum and the resulting controlled system exhibits a symmetry that has to be taken into account in the numerical approach in order to avoid numerical ill-conditioning.

\section{Mathematical formulation for optimal control of rigid bodies}\label{sec:MFOCRB}

The dynamics of rigid bodies exhibit important geometric features. The configuration of a rigid body can be described by the position vector of its center of mass in the Euclidean space $\Re^3$ and by the attitude of the rigid body represented by a rotation matrix in the special orthogonal group $\SO=\{ R\in \Re^{3\times 3}\,|\, R^TR=I,\det{R}=1\}$. Thus, the general motion of a rigid body is described by the special Euclidean group $\SE=\SO\textcircled{s}\Re^3$.  The configuration manifold for the class of multiple rigid bodies can be represented as a product involving $\Re^3,\SO$, and $\SE$. Therefore, the configuration manifold of  rigid bodies is a Lie group. Furthermore, the dynamics of rigid bodies, viewed as Lagrangian or Hamiltonian systems, are characterized by symplectic, momentum and energy preserving properties. These geometric features determine the qualitative behavior of the rigid body dynamics.

In this paper, we study optimal control problems for rigid bodies while carefully considering the geometric features of the dynamics in both the analysis and numerical computations. In particular, discrete-time dynamics of rigid bodies are developed that evolve on the configuration manifold according to a discrete version of Hamilton's principle. The resulting geometric numerical integrator, referred to as a Lie group variational integrator, preserves geometric features of the dynamics and guarantees evolution on the configuration manifold. Based on the discrete-time rigid bodies dynamics, a discrete-time optimal control problem for rigid bodies is formulated. Standard numerical optimization algorithms can then be applied to solve this discrete-time optimal control problem.  

Thus, our approach to discrete-time optimal control is characterized by discretizing the continuous-time optimal control problem at the problem formulation stage using Lie group variational integrators. This is in contrast to traditional techniques wherein discretization only arises at the last stage when numerically solving the continuous-time optimality conditions. Since the geometric properties of the dynamics of rigid bodies are preserved by using a Lie group variational integrator, this optimal control approach yields geometrically-exact optimal control inputs and accurate trajectories that are efficiently computed~\cite{Blo.ICoDaC1998,Jun.IC2005,Lee.JoDaCS2007,Lee.PotICoDaC2006}. 

In this section, we first describe the fundamental procedure to develop a Lie group variational integrator and its computational properties. Then, a discrete-time optimal control problem is formulated using the Lie group variational integrator, and computational approaches are presented to solve it numerically.

\subsection{Lie group variational integrator}\label{subsec:lgvi}
Geometric numerical integrators are numerical integration algorithms that preserve  features of the continuous-time dynamics such as invariants, symplecticity, and the configuration manifold~\cite{HaLuWa2006}. The geometrically exact properties of the discrete-time flow generate improved qualitative behavior.  In this paper, we view a Lie group variational integrator as an intrinsically discrete-time dynamical system.

Numerical integration methods that preserve the simplecticity of a Hamiltonian system have been studied extensively~\cite{San.AN92,LeRe2004}. One traditional approach is to carefully choose the  coefficients of a Runge-Kutta method to satisfy a simplecticity criterion and order conditions in order to obtain a symplectic Runge-Kutta method. However, it can be difficult to construct such integrators, and it is not guaranteed that other invariants of the system, such as momentum maps, are preserved. Alternatively, variational integrators are constructed by discretizing Hamilton's principle, rather than discretizing the continuous Euler-Lagrange equation~\cite{MosVes.CMP91,MaWe2001}. The resulting integrators have the desirable property that they are symplectic and momentum preserving, and they exhibit good energy behavior for exponentially long times. Lie group methods are numerical integrators that preserve the Lie group structure of the configuration manifold~\cite{IsMuNoZa2000}. Recently, these two approaches have been unified to obtain Lie group variational integrators that preserve the geometric properties of the dynamics as well as the Lie group structure of the configuration manifold without the use of local charts, reprojections, or constraints ~\cite{MaPeSh1999,Leo.Phd04,CMA07,CMDA07}.

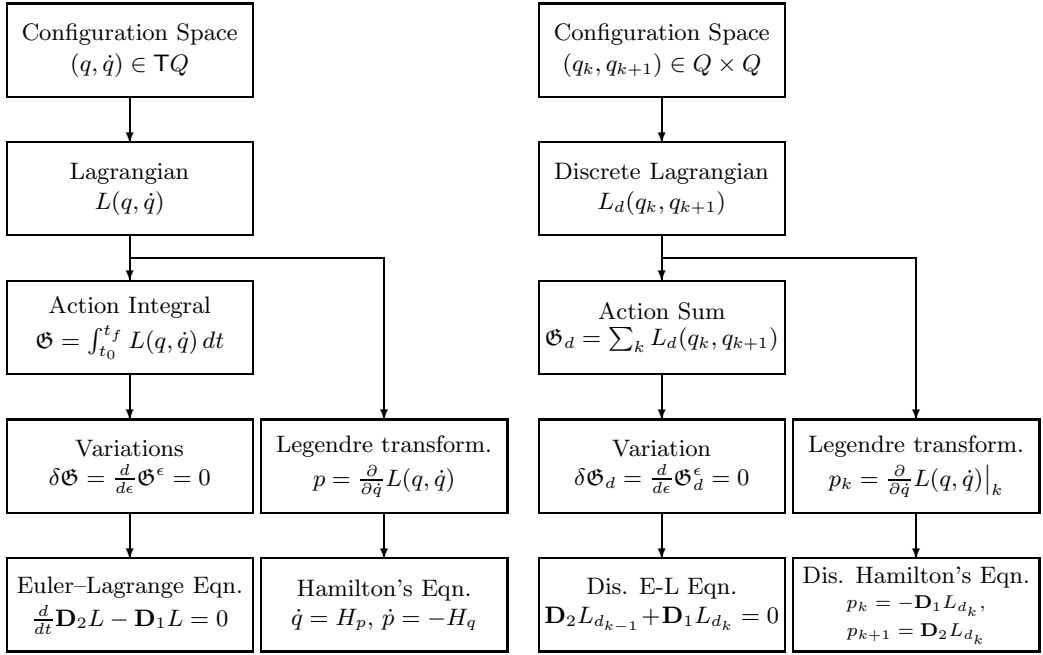
\begin{figure}
\setlength{\unitlength}{1.75em}\centering\small
\begin{picture}(22.3,14)(0,-14)
\put(0.0,-2.0){\framebox(5.3,2.0)[c]
{\shortstack[c]{Configuration Space\\$(q,\dot{q})\in \T Q$}}}
\put(2.65,-2.0){\vector(0,-1){1.0000}}
\put(0.0,-5.0){\framebox(5.3,2.0)[c]
{\shortstack[c]{Lagrangian\\$L(q,\dot{q})$}}} \put(2.65,-5.0){\vector(0,-1){1.0000}}
\put(2.65,-5.5){\line(1,0){5.5}}
\put(8.15,-5.5){\vector(0,-1){3.5}}
\put(0.0,-8.0){\framebox(5.3,2.0)[c]
{\shortstack[c]{Action Integral\\$\mathfrak{G}=\int_{t_0}^{t_f} L(q,\dot{q})\, dt$}}}
\put(2.65,-8.0){\vector(0,-1){1.0000}}
\put(0.0,-11.0){\framebox(5.3,2.0)[c]
{\shortstack[c]{Variations\\$\delta\mathfrak{G}=\frac{d}{d\epsilon}\mathfrak{G}^\epsilon=0$}}}
\put(2.65,-11.0){\vector(0,-1){1.0000}}
\put(0.0,-14.0){\framebox(5.3,2.0)[c]
{\shortstack[c]{Euler--Lagrange Eqn.\\$\frac{d}{dt}\D_2 L -\D_1 L =0$}}}
\put(5.5,-11.0){\framebox(5.3,2.0)[c]
{\shortstack[c]{Legendre transform.\\$p=\deriv{}{\dot q}L(q,\dot{q})$}}}
\put(8.15,-11.0){\vector(0,-1){1.0000}}
\put(5.5,-14.0){\framebox(5.3,2.0)[c]
{\shortstack[c]{Hamilton's Eqn.\\$\dot{q}=H_p,\,\dot{p}=-H_q$}}}
\put(11.5,-2.0){\framebox(5.3,2.0)[c]
{\shortstack[c]{Configuration Space\\$(q_k, q_{k+1})\in Q\times Q$}}}
\put(14.15,-2.0){\vector(0,-1){1.0000}}
\put(11.5,-5.0){\framebox(5.3,2.0)[c]
{\shortstack[c]{Discrete Lagrangian\\$L_{d}(q_k,q_{k+1})$}}}
\put(14.15,-5.0){\vector(0,-1){1.0000}}
%
\put(14.15,-5.5){\line(1,0){5.5}}
\put(19.65,-5.5){\vector(0,-1){3.5}}
\put(11.5,-8.0){\framebox(5.3,2.0)[c]
{\shortstack[c]{Action Sum\\$\mathfrak{G}_d=\sum_{k} L_{d}(q_k,q_{k+1})$}}}
\put(14.15,-8.0){\vector(0,-1){1.0000}}
\put(11.5,-11.0){\framebox(5.3,2.0)[c]
{\shortstack[c]{Variation\\$\delta\mathfrak{G}_d=\frac{d}{d\epsilon}\mathfrak{G}_d^\epsilon=0$}}}
\put(14.15,-11.0){\vector(0,-1){1.0000}}
\put(11.5,-14.0){\framebox(5.3,2.0)[c]
{\shortstack[c]{Dis. E-L Eqn.
\\$\D_2 L_{d_{k-1}}\!+\!\D_1 L_{d_k}=0$}}}
\put(17,-11.0){\framebox(5.3,2.0)[c]
{\shortstack[c]{Legendre transform.\\$p_k=\deriv{}{\dot q}L(q,\dot{q})\big|_k$}}}
\put(19.65,-11.0){\vector(0,-1){1.0000}}
\put(17,-14.0){\framebox(5.3,2.0)[c]
{\shortstack[c]{Dis. Hamilton's Eqn.
\\\scriptsize{$p_k=-\D_1 L_{d_k},$}\\
\scriptsize{$p_{k+1}=\D_2 L_{d_k}$}}}}
\end{picture}
\caption{Procedures to derive continuous-time and discrete-time equations of motion}\label{fig:el}
\end{figure}

We now summarize the derivation of  a Lie group variational integrator. In Lagrangian mechanics, the equations of motion are derived by finding the path that extremizes the action integral, which is the integral of the Lagrangian over time. The Legendre transformation provides an alternative description that leads to Hamilton's equations. Discrete-time Lagrangian and Hamiltonian mechanics, referred to as variational integrators, have been developed by reformulating Lagrangian and Hamiltonian mechanics in a discrete-time setting~\cite{MaWe2001}.

Discrete-time mechanics has a parallel structure with the mechanics described in continuous-time, as summarized in \reffig{el}. The phase variables of the continuous-time Lagrangian are replaced by two copies of the discrete-time configuration variables and a discrete-time Lagrangian that approximates a segment of the action integral is chosen.  An action sum is defined using the discrete-time Lagrangian such that it approximates the action integral. This is the only approximation made in the development of discrete-time mechanics.  Discrete-time Euler-Lagrange equations are obtained by setting the variation of the action sum to zero. The discrete-time Legendre transformation yields the equivalent of Hamilton's equations. Lie group variational integrators are developed to preserve the structure of the Lie group configurations as well as the geometric properties of the continuous-time dynamics. The basic idea for all Lie group methods is to express the update map for group elements in the configuration manifold in terms of the group operation, so that the group structure is preserved automatically without need of parameterizations, constraints, or reprojections.

More explicitly, consider a mechanical system whose configuration manifold is a Lie group $G$ and is described by a Lagrangian $L:\T G\rightarrow \Re$. The discrete update for the configuration is chosen as
\begin{align}
    g_{k+1} = g_k f_k,\label{eqn:g1}
\end{align}
where $g_{k},g_{k+1}\in G$ are configuration variables, and the subscript $k$ denotes the value of a variable at the time $t=kh$ for a fixed timestep $h\in\Re$. The discrete-time update map is represented by a right group action of $f_k \in G$ on $g_k$.  Since the group element is updated by a group action, the group structure is preserved.

The expression for the flow map in discrete-time is obtained from the discrete variational principle on a Lie group, as presented in \reffig{el}. A discrete Lagrangian $L_d:G\times G\rightarrow \Re$ approximates the integral of the Lagrangian over a time step along the solution of the Euler-Lagrange equation
\begin{align}
L_d(g_k,f_k) & \approx \int_{kh}^{(k+1)h} L(g(t),\dot g(t)) \,dt,
\end{align}
where a curve $g(t):[kh,(k+1)h]\rightarrow G$ satisfies the Euler--Lagrange equation in the time interval $[k,(k+1)h]$ with boundary conditions $g(kh)=g_k$ and $g((k+1)h)=g_k f_k=g_{k+1}$. Analogous to the action integral, the action sum is defined as
\begin{align}
    \mathfrak{G}_d = \sum_{k=0}^{N-1} L_d(g_k,f_k).
\end{align}
The discrete Lagrange-d'Alembert principle, which is a modification of Hamilton's principle to include the effect of control inputs, states that the sum of the variation of the action sum and the virtual work done by the control inputs is zero. But, the infinitesimal variation of a Lie group element must be carefully expressed to respect the structure of the Lie group. For example, it can be expressed in terms of the exponential map $\exp:\mathfrak{g}\rightarrow G$ as
\begin{align}
    \delta g = \frac{d}{d\epsilon}\bigg|_{\epsilon=0} g \exp {\epsilon\eta} = g\eta,\label{eqn:delg}
\end{align}
for a Lie algebra element $\eta\in\mathfrak{g}$. From the discrete Lagrange-d'Alembert principle, we obtain
\begin{align}
    \delta \mathfrak{G}_d +\sum_{k=0}^{N-1}\bracket{ u^+_k \cdot \eta_{k+1} + u^-_k \cdot \eta_k}=0
\end{align}
for any $\delta g_k$, and for given discrete Lagrangian forces $u^+_{d_k},u^-_{d_k}\in \mathfrak{g}^*$.   This yields the generalized discrete Euler--Poincar\'{e} equation
\begin{multline}
    \T_e^* \Le_{f_{k}} \cdot \D_2 L_d (g_{k},f_{k}) - \mathrm{Ad}^*_{f_{k}}\cdot(\T_e^* \Le_{f_{k+1}}\cdot \D_2 L_d (g_{k+1},f_{k+1}))\\
     + \T_e^* \Le_{g_{k+1}}\cdot \D_1 L_d (g_{k+1},f_{k+1})+u^+_{d_{k-1}}+u^-_{d_k}=0.\label{eqn:DEP}
\end{multline}
Here $\Le_f:G\rightarrow G$ denotes the left translation map given by $\Le_f g = fg$ for $f,g\in G$, $\T_g \Le_f:\T_g G\rightarrow \T_{fg} G$ is the tangent map for the left translation and $\mathrm{Ad}_g:\mathfrak{g}\rightarrow\mathfrak{g}$ is the adjoint map. A dual map is denoted by a superscript ${}^*$ (see \cite{Lee.2008} for detailed definitions and developments).

This approach has been applied to the rotation group $\mathrm{SO}(3)$ and to the special Euclidean group $\mathrm{SE}(3)$ for dynamics of rigid bodies in~\cite{CMA07,CCA05,CMDA07} and the generalization to abstract Lie groups are summarized here, thereby generating a unified geometric integrator for the class of multiple generalized rigid bodies whose configuration manifold can be expressed as a Lie group, which includes products involving $\Re^3$, $\SO$, and $\SE$ as special cases.

\subsection{Discrete-time optimal control}
Optimal control problems involve finding a control input such that a certain optimality objective is achieved under prescribed constraints. Here, the control inputs are parameterized by their values at each discrete time step, and the discrete-time equations of motion, including the control inputs, are obtained from \refeqn{DEP}. Any standard numerical algorithm for constrained optimization can be applied to this discrete-time system.

An indirect approach to solving a discrete-time optimal control problem is based on solving discrete-time necessary conditions for optimality.    The resulting two-point boundary value problem can be solved by using standard numerical root finding techniques; one such approach is the shooting method that iterates on initial values of the multipliers.    Alternatively, a direct approach formulates the discrete-time optimal control problem as a nonlinear programming problem, which is solved using standard numerical optimization algorithms such as a sequential quadratic programming algorithm; one such approach is the DMOC (Discrete Mechanics and Optimal Control) approach~\cite{Jun.IC2005}.

Explicit time-discretization prior to numerical optimization has significant computational advantages. As discussed in the previous section, the discrete-time dynamics are faithful representations of the continuous-time dynamics, and consequently more accurate solutions to the optimal control problems are typically obtained. The external control inputs may break the Lagrangian and Hamiltonian system structure; for example, the total energy may not be conserved for a controlled mechanical system. But, the computational superiority of the discrete mechanics formulation still holds for controlled systems. In particular, it has been demonstrated in~\cite{MaWe2001} that the discrete-time dynamics derived from the discrete Lagrange-d'Alembert principle accurately computes the energy dissipation rate of controlled systems. For example, this feature is extremely important in accurately computing optimal trajectories for spacecraft orbit and attitude maneuvers for which the control authority is low and the maneuver time is large.

The proposed discrete-time optimal control formulation provides a framework for accurate computations. In most indirect optimal control approaches, the optimal solutions are sensitive to small variations in the initial values of the multipliers. This may cause difficulties, such as numerical ill-conditioning, in solving the necessary conditions for optimality expressed as a two-point boundary value problem. Numerically computed sensitivity derivatives, using Lie group variational integrators, do not exhibit numerical dissipation, which typically arises in conventional numerical integration schemes. Thus, the proposed approach leads to numerical robustness and efficient numerical computations.  This indirect computational approach exhibits the quadratic convergence rate that is typical of Newton methods when it is applied to an optimal attitude control problem~\cite{ACC07.opt}; the error in satisfaction of the optimality condition converges to machine precision superlinearly. For the direct optimal control approach, the optimal control inputs can be parametrized using fewer degrees of freedom, thereby reducing the computational overhead.

Several optimal control problems involving rigid bodies have been previously studied by the authors. Minimum-fuel and time-optimal control of spacecraft large-angle attitude maneuvers are studied in~\cite{ACC06,HusMelSan.CDC06,Lee.JoDaCS2007,ACC08}.  The optimal orbit transfer of a dumbbell spacecraft, wherein the rotational attitude dynamics are non-trivially coupled to the translational dynamics, is studied in~\cite{Lee.PotICoDaC2006}.  An underactuated optimal control problem for the attitude maneuver of a 3D pendulum is studied in~\cite{ACC07.opt}.   An optimal formation reconfiguration of multiple rigid body spacecraft is studied in~\cite{CDC07.opt}. An optimal control problem for a dynamic system evolving on an abstract Lie group is developed in~\cite{Lee.2008}, thereby generating a unified approach for optimal control problems of multiple rigid bodies.

In this paper, we summarize results for two optimal control problems for a single rigid body in Section~\ref{sec:optsrbdy} and results for two optimal control problems for multiple rigid bodies in Section~\ref{sec:optmrbdy}. Each of these optimal control problems treats complex dynamics of a single or multiple rigid bodies, demonstrating the value of the proposed geometric optimal control approach.

\section{Optimal control problems for a single rigid body}\label{sec:optsrbdy}

\subsection{Optimal maneuver of a dumbbell spacecraft on $\SE$}\label{subsec:optse}
We develop an optimal 3D translational and rotational maneuver of a rigid dumbbell spacecraft in orbit about a large central body. The dumbbell spacecraft is composed of two spheres connected by a massless rod. An interesting feature of the dumbbell spacecraft is that there is coupling between its translational dynamics and its rotational dynamics due to the presence of both gravity forces and gravity moments that act on the dumbbell spacecraft.

The configuration manifold is the special Euclidean group $\SE=\SO\,\textcircled{s}\,\Re^3$. For $(R,x)\in\SE$, the linear transformation from the body-fixed frame to the inertial frame is denoted by the rotation matrix $R\in\SO$, and the position of the mass center in the inertial frame is denoted by a vector $x\in\Re^3$. The vectors $\Omega,v\in\Re^3$ are the angular velocity in the body-fixed frame, and the translational velocity in the inertial frame, respectively. Let $m\in\Re$ and $J\in\Re^{3\times 3}$ be the mass and the moment of inertia matrix of a rigid body. We assume that external control force $u^f\in\Re^3$ and control moment $u^m\in\Re^3$ act on the dumbbell spacecraft. Control inputs are parameterized by their values at each time step.

Define a $f_k=(F_k,Y_k)\in\SE$ such that $g_{k+1}=(R_{k+1},x_{k+1})$ is equal to $g_kf_k$, i.e. $(R_{k+1},x_{k+1})=(R_k,x_k)\circ(F_k,Y_k)= (R_kF_k, x_k+R_kY_k)$. The rotation matrix $F_k$ represent the relative update of the attitude between integration steps. The gravitational potential is denoted by $U:\SE\rightarrow\Re$. We choose the following discrete Lagrangian
\begin{align}
L_d(R_k,x_k,F_k,Y_k) = \frac{1}{2h}mY_k^TY_k+\frac{1}{h}\tr{(I-F_k)J_d}-hU(R_{k}F_k,x_{k}+R_kY_k),
\end{align}
where $J_d\in\Re^{3\times 3}$ is a non-standard moment of inertia matrix defined as $J_d=\frac{1}{2}\mathrm{tr}[J]I_{3\times 3}-J$. Substituting this discrete Lagrangian into \refeqn{DEP}, we obtain the following discrete equations of motion (see \cite{Lee.2008} for detailed development).
\begin{gather}
h \widehat{J\Omega_{k}}=F_{k}J_{d}-J_{d}F_{k}^T,\label{eqn:findF}\\
R_{k+1}=R_{k}F_{k},\label{eqn:Rkp}\\
x_{k+1} = x_{k} + hv_k,\label{eqn:xkp}\\
J\Omega_{{k+1}}=F_{k}^T J\Omega_{k}+h (M_{{k+1}}+u^m_{k+1})
\label{eqn:Omegakp},\\
m v_{k+1}=m v_k-h\deriv{U_{k+1}}{x_{k+1}}+hu^f_{k+1},\label{eqn:vkp}
\end{gather}
where the hat map $\hat\cdot$ is an isomorphism from $\Re^3$ to $3\times 3$ skew-symmetric matrices $\so$, defined such that $\hat x y =x\times y$ for any $x,y\in\Re^3$. The moment $M\in\Re^3$ due to the potential is given by,
\begin{align}
\hat M & = \deriv{U}{R}^TR-R^T\deriv{U}{R},
\end{align}
where the matrix $\deriv{U}{R}\in\Re^{3\times 3}$ is defined by $[\deriv{U}{R}]_{ij}=\deriv{U}{[R]_{ij}}$ for $i,j\in\braces{1,2,3}$, and the $i,j$-th element of a matrix is denoted by $[\cdot]_{ij}$.

For a given $(R_k,x_k,\Omega_k,v_k)$, we solve the implicit equation \refeqn{findF} to find $F_k\in\SO$. Then, the configuration at the next step $(R_{k+1},x_{k+1})$ is obtained from \refeqn{Rkp} and \refeqn{xkp}. Using the computed moment $M_{k+1}$ and force $-\deriv{U_{k+1}}{x_{k+1}}$, velocities $\Omega_{k+1},v_{k+1}$  are obtained from \refeqn{Omegakp} and \refeqn{vkp}. This defines a discrete flow map, $(R_k,x_k,\Omega_k,v_k)\mapsto(R_{k+1},x_{k+1},\Omega_{k+1},v_{k+1})$, and this process can be repeated.

Since this Lie group variational integrator is obtained by discretizing Hamilton's principle, it is symplectic and preserves the momentum map associated with the symmetry of the Lagrangian. In the absence of external forces and moments, the total energy oscillates around its initial value with small bounds on a comparatively short timescale, but there is no tendency for the mean of the oscillation in the total energy to drift (increase or decrease) from the initial value for exponentially long times.

The discrete flow map also preserves the group structure. By using the given computational approach, the matrix $F_k$, representing the change in relative attitude change over a time step, is guaranteed to be a rotation matrix. The rotation matrix $R_{k+1}$ is obtained by the group operation in \refeqn{Rkp}, so that it evolves on $\SO$. Therefore, the orthogonal structure of the rotation matrices is preserved, and the attitude of each rigid body is determined accurately and globally.

This geometrically exact numerical integration method yields a highly efficient computational algorithm. The self-adjoint discrete Lagrangian used to derive this Lie group variational integrator guarantees that this integrator has second-order accuracy, while requiring only one function evaluation per integration step. Higher-order methods can be easily constructed using a composition method~\cite{HaLuWa2006}.

An implicit equation \refeqn{findF} must be solved at each time step to determine the attitude update. However the computational effort to solve each implicit equation is negligible; the relative attitude update is expressed at the Lie algebra level isomorphic to $\Re^3$, and the corresponding Newton iteration converges to machine precision within two or three iterations. This method could be considered \textit{almost explicit} when the computational cost is compared with explicit integrators with the same order of accuracy~\cite{CMA07}. 

\subsubsection*{Optimal control problem}
The objective is to transfer the spacecraft from a given initial condition $(R_0,x_0,\Omega_0,v_0)$ to a desired terminal condition $(R^f,x^f,\Omega^f,v^f)$ during a fixed maneuver time $Nh$, while minimizing the square of the $l_2$ norm of the control inputs.
\begin{gather}
\min_{u_{k+1}} \braces{\mathcal{J}=\sum_{k=0}^{N-1}
\frac{h}{2}(u^f_{k+1})^TW_fu^f_{k+1}+
\frac{h}{2}(u^m_{k+1})^TW_mu^m_{k+1}},
\end{gather}
where $W_f,W_m\in\Re^{3\times 3}$ are symmetric positive-definite matrices.

\subsubsection*{Necessary conditions for optimality}
An indirect optimization method is used to determine the optimal solution, based on necessary conditions for optimality derived using variational arguments; the optimal control is characterized as a solution of a two-point boundary value problem. The augmented cost function to be minimized is
\begin{align}
\mathcal{J}_a =
\sum_{k=0}^{N-1}&\frac{h}{2}(u^f_{k+1})^TW^fu^f_{k+1}+
\frac{h}{2}(u^m_{k+1})^TW^mu^m_{k+1}\nonumber\\
& +\lambda_k^{1,T}\braces{-x_{k+1}+x_k+h v_k} +\lambda_k^{2,T}\braces{-m v_{k+1} + mv_k-h\deriv{U_{k+1}}{x_{k+1}}+hu^f_{k+1}}\nonumber\\
& +\lambda_k^{3,T}\parenth{\mathrm{logm}(F_k-R_{k}^TR_{k+1})}^{\vee} +\lambda_k^{4,T}\braces{-J\Omega_{k+1} + F_k^T J\Omega_k +
h\parenth{M_{k+1}+u_{k+1}^m}},
\end{align}
where $\lambda_k^1,\lambda_k^2,\lambda_k^3,\lambda_k^4\in \Re^3$ are Lagrange multipliers. The matrix logarithm is denoted by $\mathrm{logm}:\SO\rightarrow\so$ and the vee map $\vee:\so\rightarrow\Re^3$ is the inverse of the hat map. The logarithmic form of \refeqn{Rkp} is used, and the constraint \refeqn{findF} is implicitly imposed using constrained variations. Using  similar expressions for the variations given in \refeqn{delg}, the infinitesimal variation of the cost can be written as
\begin{align}
\delta\mathcal{J}_a & = \sum_{k=1}^{N-1} h\delta
u_{k}^{f,T}\braces{W_fu^f_{k}+\lambda_{k-1}^2}
+h\delta u_{k}^{m,T}\braces{W_mu^m_{k}+\lambda_{k-1}^4}
+z_k^T\braces{-\lambda_{k-1}+A_k^T \lambda_k},
\end{align}
where $\lambda_k=[\lambda_k^1;\lambda_k^2;\lambda_k^3;\lambda_k^4]\in\Re^{12}$ is the vector of Lagrange multipliers, and $z_k\in\Re^{12}$ represents the infinitesimal variation of $(R_k,x_k,\Omega_k,v_k)$, given by $z_k=[\mathrm{logm}(R_k^T \delta R_k)^\vee;\delta x_k,\delta\Omega_k,\delta v_k]$. The matrix $A_k\in\Re^{12\times 12}$ is expressed in terms of $(R_k,x_k,\Omega_k,v_k)$~\cite{Lee.PotICoDaC2006}. Thus, necessary conditions for optimality are given by
\begin{align}
u^f_{k+1} &= -W_{f}^{-1}\lambda_{k}^2,\label{eqn:ufkp}\\
u^m_{k+1} &= -W_{m}^{-1}\lambda_{k}^4,\label{eqn:umkp}\\
\lambda_{k} &= A_{k+1}^T \lambda_{k+1}\label{eqn:updatelam}
\end{align}
together with the discrete equations of motion and the boundary conditions.

\subsubsection*{Computational approach}
Necessary conditions for optimality are expressed in terms of a two-point boundary problem. This problem is to find the optimal discrete flow, multipliers, and control inputs that simultaneously satisfies the equations of motion, optimality conditions, multiplier equations, and boundary conditions. We use a neighboring extremal method~\cite{Bry.BK75}, and choose a nominal solution satisfying all of the necessary conditions except the boundary conditions. The unspecified initial multiplier is updated by successive linearization so as to satisfy the specified terminal boundary conditions in the limit. This is also referred to as a shooting method. The main advantage of the neighboring extremal method is that the number of iteration variables is small.

The difficulty is that the extremal solutions are sensitive to small changes in the unspecified initial multiplier values. The nonlinearities also make it hard to construct an accurate estimate of sensitivity, thereby resulting in numerical ill-conditioning. Therefore, it is important to compute the sensitivities accurately in the neighboring extremal method. Here, the optimality conditions \refeqn{ufkp} and \refeqn{umkp} are substituted into the equations of motion and the multiplier equations, which are linearized to obtain 
\begin{align*}
    \begin{bmatrix} z_N \\ \delta \lambda_N \end{bmatrix}
    =\begin{bmatrix} \Psi^{11} & \Psi^{12} \\ \Psi^{21} & \Psi^{22}\end{bmatrix}
    \begin{bmatrix} z_0 \\ \delta \lambda_0 \end{bmatrix},
\end{align*}
where $\Psi^{ij}\in\Re^{6\times 6}$ for $i,j\in\{1,2\}$ represents a computable linear operator. For the given two-point boundary value problem, $z_0=0$ since the initial condition is fixed. The terminal multipliers are free. Thus, we obtain
\begin{align*}
    z_N = \Psi^{12} \delta\lambda_0.
\end{align*}
The linear operator $\Psi^{12}$ represents the sensitivity of the specified terminal boundary conditions with respect to the unspecified initial multiplier. Using this sensitivity, a guess of the unspecified initial multipliers is iterated to satisfy the specified terminal conditions in the limit. Any type of Newton iteration can be applied. We use a line search with backtracking algorithm, referred to as the Newton-Armijo iteration~\cite{Kel.BK95}.

\subsubsection*{Numerical example}

We study a maneuver of a rigid spacecraft under a central gravity field. We assume that the mass of the spacecraft is negligible compared to the mass of a central body, and we consider a fixed frame attached to the central body as an inertial frame. The resulting model is a Restricted Full Two Body Problem (RF2BP)~\cite{Sch.CMDA02}.

The spacecraft is modeled as a dumbbell, which consists of two equally massive spheres and a massless rod. The gravitational potential is given by
\begin{align}
U(R,x)=-\frac{GMm}{2} \sum_{q=1}^2 \frac{1}{\norm{x+R\rho^q}},
\end{align}
where $G\in\Re$ is the gravitational constant, $M,m\in\Re$ are the mass of the central body, and the mass of the dumbbell, respectively. The vector $\rho^q\in\Re^3$ is the position of the $q$th sphere from the mass center of the dumbbell expressed in the body fixed frame ($q\in \braces{1,2}$). The mass, length, and time dimensions are normalized by the mass of the dumbbell, the radius of a reference circular orbit, and its orbital period.

Initially, the spacecraft is on a circular orbit. The desired maneuver is to increase the orbital inclination by $60^\circ$. We explicitly consider the coupling effect between the orbital motion and the rotational attitude maneuver of the spacecraft. The maneuver time is chosen to be a quarter of the orbital period of the initial circular orbit. The boundary conditions are as follows,
\begin{alignat*}{2}
x_0&=[1,0,0],&\quad x^f&=[-0.3536,0.3536,0.8660],\\
R_{0}&=\begin{bmatrix}
     0&     -1&     0\\
     1&     0&     0\\
     0&     0&     1\end{bmatrix},&
\quad R^f&=\begin{bmatrix}
     -0.7071&     0.3535&     0.6123\\
     -0.7071&     -0.3535&    -0.6123\\
     0&     -0.8660&     0.5\end{bmatrix},\\
\dot x_0&=[0,0.9835,0],&\quad \dot x^f&=[-0.6954,-0.6954,0],\\
\Omega_0&=[0,0,0.9835],&\quad\Omega^f&=[0,0,0.9835].
\end{alignat*}

\begin{figure}[h]
\centerline{
    \subfigure[Spacecraft maneuver]{
        \includegraphics[width=0.40\textwidth]{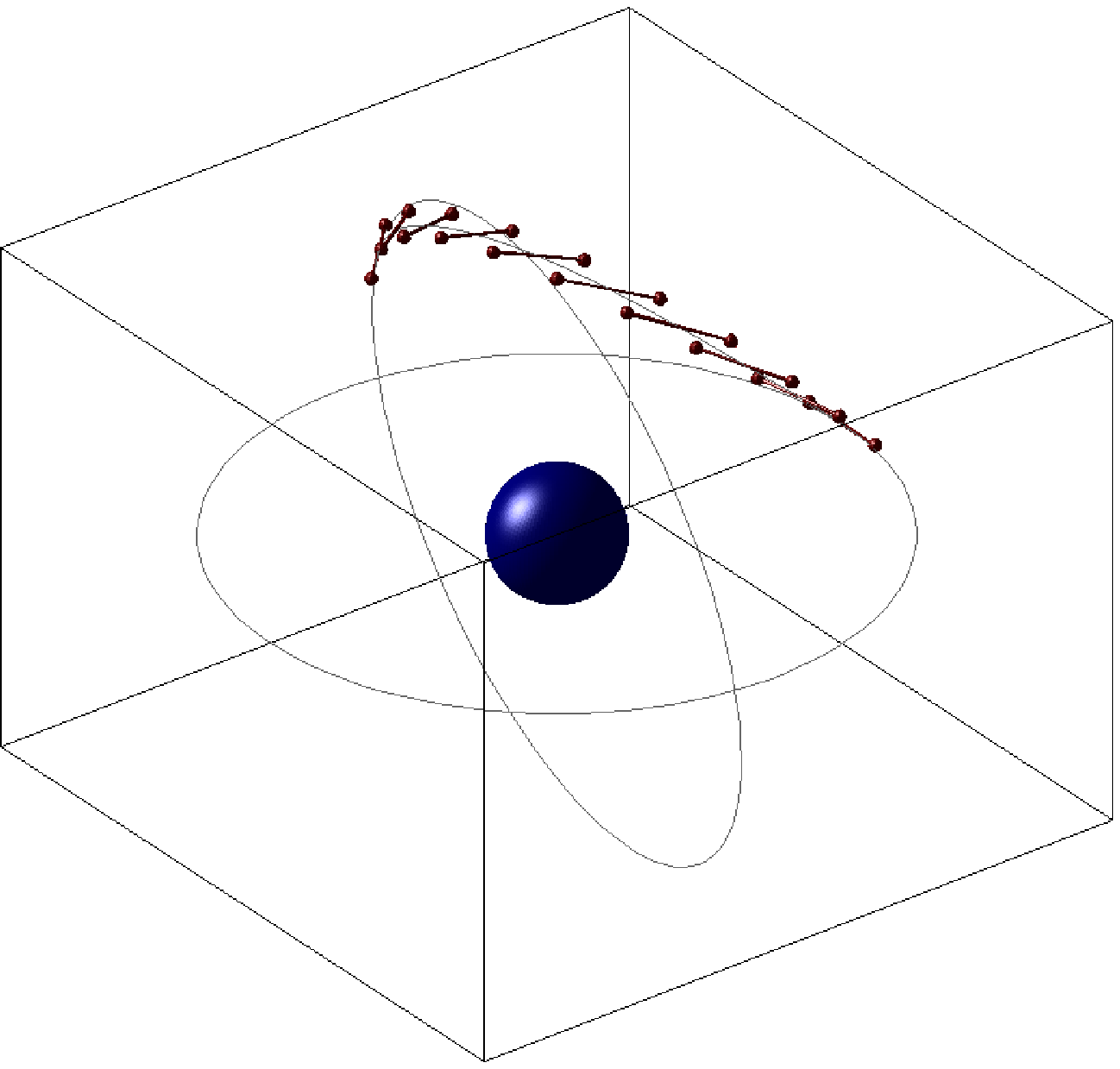}}
    \hfill
    \subfigure[Convergence rate]{
        \includegraphics[width=0.45\textwidth]{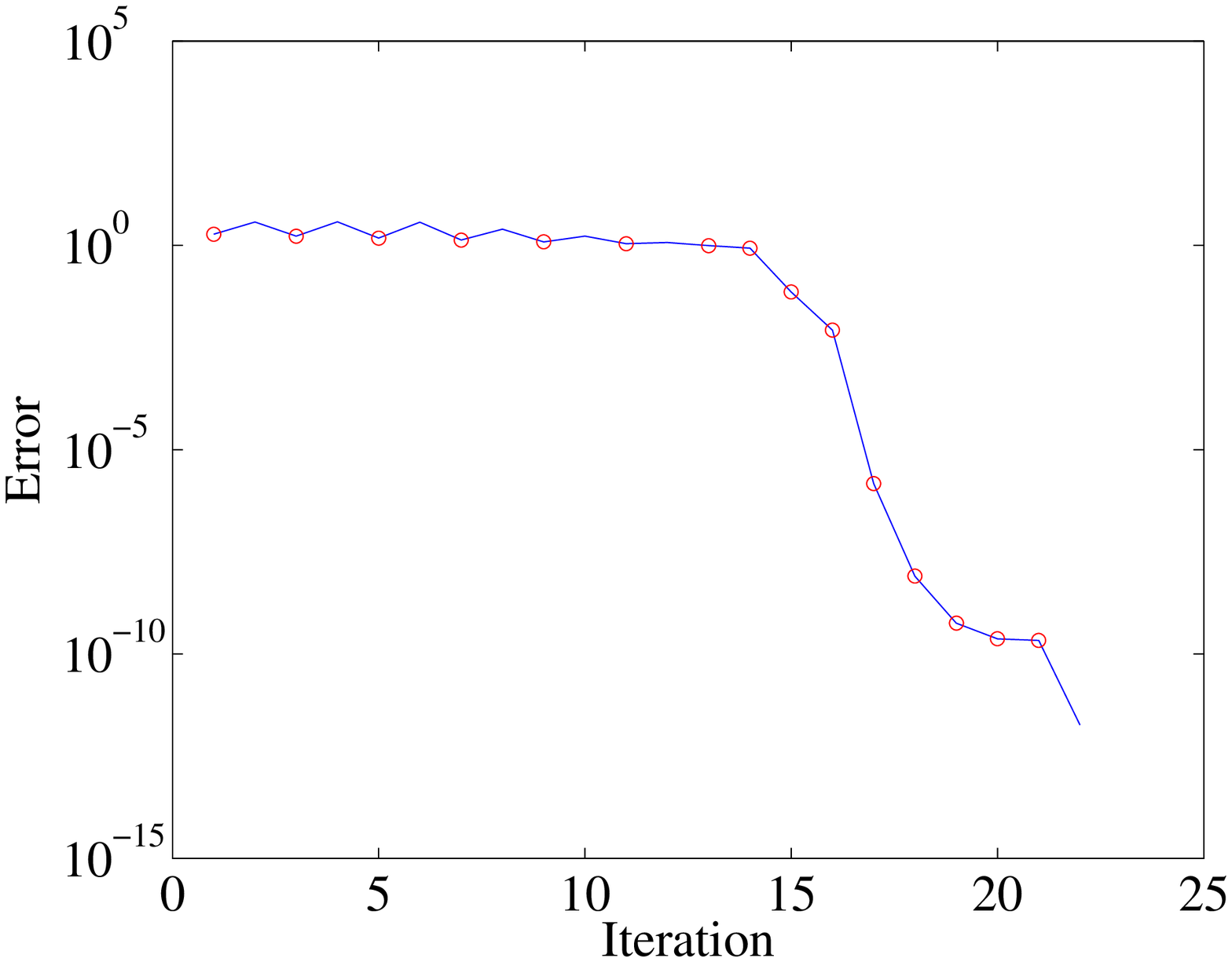}\label{fig:OptRB_err}}
}
\centerline{
    \subfigure[Control force $u^f$]{
        \includegraphics[width=0.43\textwidth]{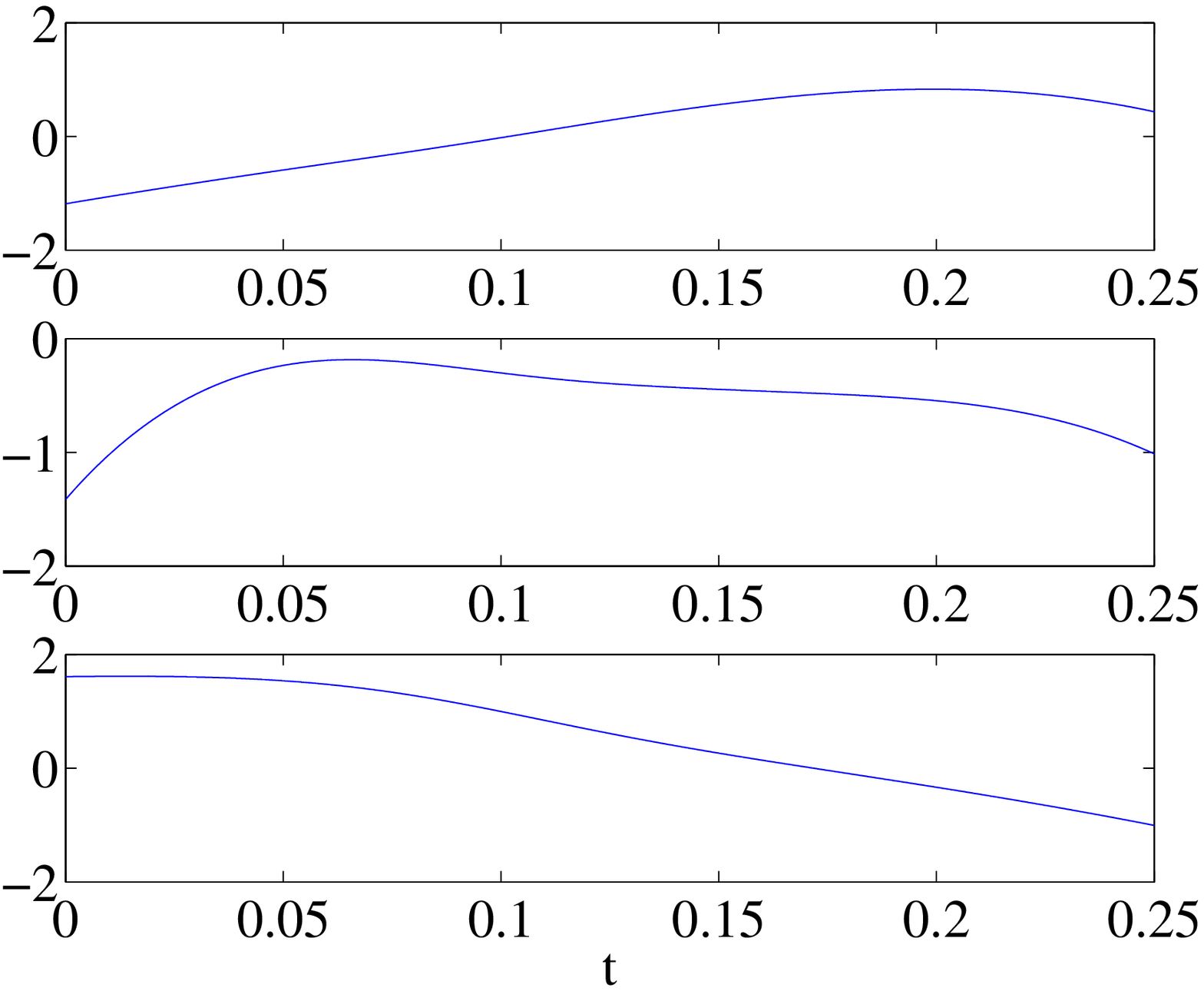}}
    \hfill
    \subfigure[Control moment $u^m$]{
        \includegraphics[width=0.455\textwidth]{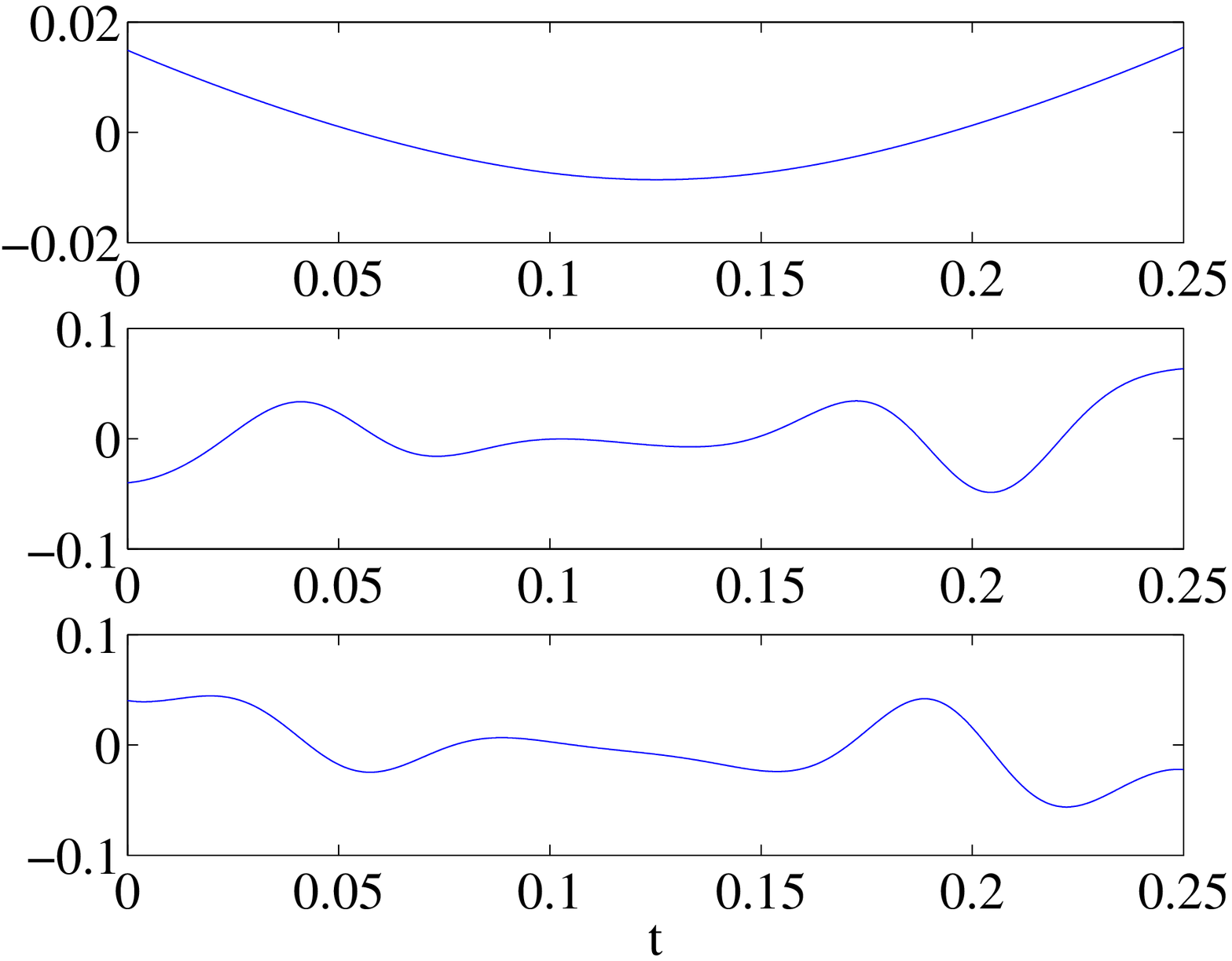}}
    }
\caption{Optimal orbit transfer of a dumbbell spacecraft}\label{fig:OptRB}
\end{figure}

\reffig{OptRB} illustrates the optimal spacecraft maneuver, convergence rate, and optimal control inputs. The optimal cost and the violation of the terminal boundary conditions are $13.03$, and $9.32\times 10^{-15}$ respectively. \reffig{OptRB_err} shows the violation of the terminal boundary conditions versus the number of iterations on a semi-logarithmic scale. Red circles denote outer iterations of the Newton-Armijo iteration where the sensitivity derivatives are computed, and inner iterations correspond to backtracking in the line search routine. The initial guess of the unspecified initial multipliers is arbitrarily chosen. The error in satisfaction of the terminal boundary condition converges quickly to machine precision after the 20th iteration. These convergence results are consistent with the quadratic convergence rates expected of Newton methods with accurately computed gradients.

The shooting method may be prone to numerical ill-conditioning, as a small change in the initial multiplier can cause highly nonlinear behavior of the terminal conditions. However, as shown in \reffig{OptRB_err}, the computational geometric optimal control approach exhibits excellent numerical convergence properties. This is because the proposed computational algorithms are geometrically exact and numerically accurate. There is no numerical dissipation introduced by the numerical algorithm, and therefore, the sensitivity derivatives are more accurately computed.

\subsection{Optimal attitude reorientation of an underactuated 3D pendulum on $\SO$ \cite{ACC07.opt}}\label{subsec:3dpend}

A 3D pendulum is a rigid body supported by a fixed frictionless pivot acting under the influence of a uniform gravitational field~\cite{SheSanNal.CDC04}. The rigid body has three rotational degrees of freedom, and the configuration manifold is $\SO$. The linear transformation from the body fixed frame and the inertial frame is denoted by $R\in\SO$, and the angular velocity represented in the body fixed frame is denoted by $\Omega\in\Re^3$. Let $e_3\in\Re^3$ be the gravity direction in the inertial frame, and $J\in\Re^{3\times 3}$ be the moment of inertia matrix of the rigid body with respect to the pivot point. The vector from the pivot point to the mass center, represented in the body fixed frame is given by $\rho\in\Re^3$.

The Lagrangian of the 3D pendulum is invariant under a rotation about the gravity direction, and therefore the 3D pendulum has a $S^1$ symmetry action. Consequently, the angular momentum about the gravity direction, represented by $e_3^T RJ\Omega$, is preserved.

We study an optimal attitude control of the 3D pendulum with symmetry. An external control moment is chosen such that it does not have any component about the gravity direction. The structure of the control moment is chosen as $R^T e_3 \times u$ for a control parameter $u\in\Re^3$. Thus, the angular momentum about the gravity direction is conserved along the controlled dynamics of the 3D pendulum. Such control inputs are physically realized by actuation mechanisms, such as point mass actuators, that change the center of mass of the 3D pendulum.

The discrete Lagrangian of the 3D pendulum is chosen to be
\begin{align*}
    L_d(R_k,F_k) = \frac{1}{h} \tr{(I-F_k)J_d} + h mg e_3^T R\rho.
\end{align*}
The resulting Lie group variational integrator, including an external control input, is given by
\begin{gather}
h \widehat{J\Omega_{k}}=F_{k}J_{d}-J_{d}F_{k}^T,\label{eqn:findF2}\\
R_{k+1}=R_{k}F_{k},\label{eqn:Rkp2}\\
J\Omega_{{k+1}}=F_{k}^T J\Omega_{k}+hM_{{k+1}}+hR^T e_3 \times u_{k+1}.\label{eqn:Omegakp2}
\end{gather}

\subsubsection*{Optimal control problem}

The objective of the optimal control problem is to transfer the 3D pendulum from a given initial condition $(R_0,\Omega_0)$ to a desired terminal condition $(R^f,\Omega^f)$ during a fixed maneuver time $Nh$, while minimizing the square of the $l_2$ norm of the control inputs.
\begin{gather}
\min_{u_{k+1}} \braces{\mathcal{J}=\sum_{k=0}^{N-1}
\frac{h}{2}(u_{k+1})^TWu_{k+1}},
\end{gather}
where $W\in\Re^{3\times 3}$ is a symmetric positive-definite matrix. In particular, we choose attitude maneuvers that can be described by rest-to-rest rotations about the unactuated gravity direction. The resulting optimal attitude maneuver exhibits the geometric phase effect~\cite{MarRat.BK99}, which in the zero group momentum case directly relates the group motion to the curvature enclosed by the trajectory in shape space.

\subsubsection*{Necessary conditions for optimality}
We solve this optimal control problem by using an indirect optimization method, where necessary conditions for optimality are derived using variational arguments, and a solution of the corresponding two-point boundary value problem provides the optimal control. The augmented cost function to be minimized is
\begin{align}
\mathcal{J}_a =
\sum_{k=0}^{N-1}& \frac{h}{2}u_{k+1}^TWu_{k+1} +\lambda_k^{1,T}\parenth{\mathrm{logm}(F_k-R_{k}^TR_{k+1})}^{\vee}\nonumber \\ & +\lambda_k^{2,T}\braces{-J\Omega_{k+1} + F_k^T J\Omega_k +
hM_{k+1}+hR_{k+1}^Te_3\times u_{k+1}},
\end{align}
where $\lambda_k^1,\lambda_k^2\in \Re^3$ are Lagrange multipliers. The infinitesimal variation can be written as
\begin{align}
\delta\mathcal{J}_a & = \sum_{k=1}^{N-1} h\delta u_{k}^{T}\braces{Wu_{k}-R_k^Te_3\times \lambda_{k-1}^2}
+z_k^T\braces{-\lambda_{k-1}+A_k^T \lambda_k},
\end{align}
where $\lambda_k=[\lambda_k^1;\lambda_k^2]\in\Re^{6}$, and $z_k\in\Re^{6}$ represents the infinitesimal variation of $(R_k,\Omega_k)$, given by $z_k=[\mathrm{logm}(R_k^T \delta R_k)^\vee;\delta\Omega_k]$. The matrix $A_k\in\Re^{6\times 6}$ can be expressed in terms of $(R_k,\Omega_k),\lambda_k$. Thus, necessary conditions for optimality are given by
\begin{align}
u_{k+1} &= W^{-1}(R_{k+1}^Te_3\times \lambda_{k}^2),\\
\lambda_{k} &= A_{k+1}^T \lambda_{k+1}
\end{align}
together with the discrete equations of motion and the boundary conditions.

\subsubsection*{Computational approach}

We apply the neighboring extremal method described in Section~\ref{subsec:optse}; the optimality condition is substituted into the equations of motion and the multiplier equation, and sensitivity derivatives of the optimal solution with respect to the initial multiplier are obtained, and the initial multiplier is iterated to satisfy the terminal boundary condition.

However, the underactuated control input, that respects the symmetry of the 3D pendulum, causes a fundamental singularity in the sensitivity derivatives, since the controlled system inherits the $S^1$ symmetry, and the cost functional is invariant under the lifted action of $S^1$. Consequently, the sensitivity derivatives vanish in the group direction. At each iteration, we need to compute inverse of a matrix of sensitivity derivatives to update the initial multiplier. However, the sensitivity matrix has a theoretical rank deficiency of one since the vertical component of the inertial angular momentum is conserved regardless of the initial multiplier variation. Therefore, this matrix inversion is numerically ill-conditioned.

We present a simple numerical scheme to avoid the numerical ill-conditioning caused by this symmetry. At each step, we decompose the matrix of sensitivity derivatives into a symmetric part and an anti-symmetric part. The symmetric part describes the sensitivity of the conserved angular momentum component due to the symmetry, and therefore it is zero and does not depend on the initial multiplier values.    An update for the initial multipliers is determined using the matrix inverse of the anti-symmetric part; this matrix inverse is not ill-conditioned. This approach removes the singularity in the sensitivity derivatives completely, and the resulting optimal control problem is no longer ill-conditioned.

\subsubsection*{Numerical example}

Properties of the 3D pendulum are chosen as,
\[ m=1\,\mathrm{kg},\quad J=\mathrm{diag}[0.13,0.28,0.17]\,\mathrm{kgm^2},\quad \text{and}\quad\rho=[0,0,0.3]\,\mathrm{m}.\]
The desired maneuver is a $180^\circ$ rotation about the vertical axis from a hanging equilibrium to another hanging equilibrium. The corresponding boundary conditions are given by
\begin{alignat*}{2}
    R_0&=I,&\quad R^f&=\mathrm{diag}[-1,-1,1],\\
    \Omega_0&=[0,0,0],&\quad \Omega^f&=[0,0,0].
\end{alignat*}
The maneuver time is $1$ second, and the time step is $h=0.001$. Since the vertical component of the angular momentum is zero, the rotation is a consequence of the geometric phase effect~\cite{MarRat.BK99}. This problem is challenging in the sense that the desired maneuvers are rotations about the gravity direction, but the control input does not directly generate any moment about the gravity direction.

\begin{figure}[h]
\centerline{
    \subfigure[3D pendulum maneuver]{
        \includegraphics[width=0.90\textwidth]{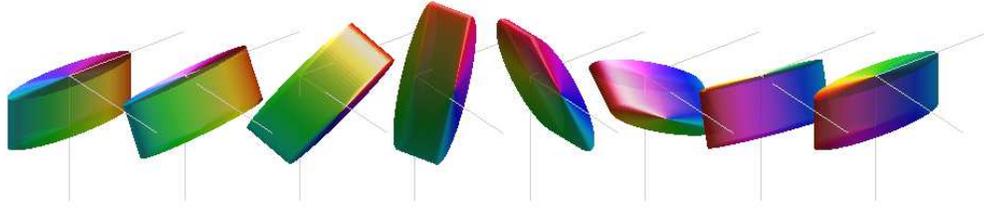}}
}
\centerline{
    \subfigure[Control moment $u$]{
        \includegraphics[width=0.45\textwidth]{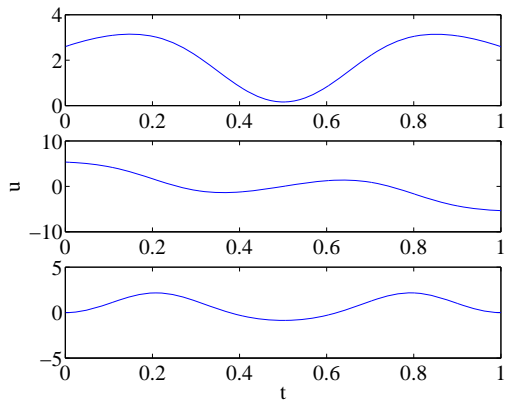}}
    \hfill
    \subfigure[Convergence rate]{
        \includegraphics[width=0.45\textwidth]{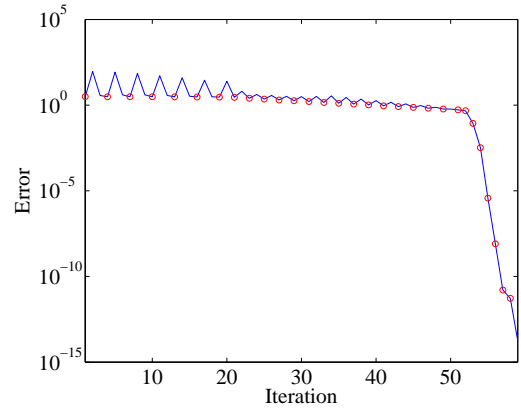}\label{fig:Opt3DP_err}}
    }
\caption{Optimal control of a 3D pendulum with symmetry}\label{fig:Opt3DP}
\end{figure}

\reffig{Opt3DP} illustrates the optimal pendulum maneuver, convergence rate, and optimal control inputs. The optimal cost and the violation of the terminal boundary conditions are $7.32$, and $4.80\times 10^{-15}$ respectively. As shown in \reffig{Opt3DP_err}, the error in satisfaction of the terminal boundary condition converges to machine precision after the 50th iteration. The condition number of the decomposed sensitivity derivative varies from $10^0$ to $10^5$. If the sensitivity derivative is not decomposed, then the condition numbers are at the level of $10^{19}$, and the numerical iterations fail to converge. This numerical example demonstrates the excellent numerical convergence properties of the computational geometric optimal control approach that is achieved by incorporating a modification that eliminates the numerical ill-conditioning introduced by the symmetry.

\section{Optimal control problems for multiple rigid bodies}\label{sec:optmrbdy}

\subsection{Optimal maneuver of a 3D pendulum on a 2D cart on $\SO\times\Re^2$}

Consider a 3D pendulum whose pivot is attached to a cart that can translate on a horizontal plane. This is a generalization of the popular planar pendulum on a cart model (see, for example,~\cite{Blo.PotICoDaC2005}), where the pendulum has three rotational degrees of freedom, and the cart moves on a two dimensional horizontal plane.

We define two frames; an inertial frame and a body fixed frame for the 3D pendulum whose origin is located at the moving pivot point. Define
\begin{center}
\begin{tabular}{lp{13cm}}
$x\in\Re$ & Displacement of the cart along the $e_1$ direction in the reference frame\\
$y\in\Re$ & Displacement of the cart along the $e_2$ direction in the reference frame\\
$R\in\SO$ & Rotation matrix from the body fixed frame to the reference frame\\
$\Omega\in\Re^3$ & Angular velocity of the pendulum represented in the body fixed frame\\
$d\in\Re^3$ & Vector from the pivot to the mass center of the pendulum represented in the body fixed frame\\
$m\in\Re$ & Mass of the pendulum\\
$M\in\Re$ & Mass of the cart
\end{tabular}
\end{center}
The configuration manifold is $\SO\times\Re^2$. We assume that external control forces $u_x,u_y\in\Re$ are applied to the cart.

The Lagrangian of the 3D pendulum on a cart is invariant under a rotation about the gravity direction. Therefore, it has a symmetry of $\S^1$ action, and the total angular momentum about the gravity direction is preserved. The external control forces acting on the cart break this symmetry, and the controlled system is not symmetric. In particular, the total angular momentum is not preserved in the controlled dynamics. Therefore, this optimal control problem should be distinguished from the optimal control of a 3D pendulum with symmetry, discussed in Section \ref{subsec:3dpend}, where a symmetry-preserving control input is chosen.

The discrete Lagrangian for the 3D pendulum on a 2D cart is
\begin{align}
L_d &(R_k,x_k,y_k,R_{k+1},x_{k+1},y_{k+1})=\frac{1}{2h}(M+m) ((x_{k+1}-x_k)^2 + (y_{k+1}-y_k)^2)\nonumber\\ &+\frac{1}{h}\tr{(I-F_k)J_d} +\frac{m}{h}(x_{k+1}-x_k) e_1^T (R_{k+1}-R_k) d + \frac{m}{h}(y_{k+1}-y_k) e_2^T (R_{k+1}-R_k) d\nonumber\\ &+\frac{h}{2}mg e_3^T R_kd+\frac{h}{2}mg e_3^T R_{k+1}d.\label{eqn:Ld}
\end{align}
From \refeqn{DEP}, the Lie group variational integrator for the 3D pendulum on a cart is given by the discrete-time equations
{\allowdisplaybreaks
\begin{align}
p_{x_k} & = \frac{1}{h}(M+m) (x_{k+1}-x_k) +\frac{m}{h}e_1 (R_{k+1}-R_k) d,\label{eqn:pxk}\\
p_{y_k} & = \frac{1}{h}(M+m) (y_{k+1}-y_k) +\frac{m}{h}e_2 (R_{k+1}-R_k) d,\label{eqn:pyk}\\
\hat p_{\Omega_k} & = \frac{1}{h} (F_kJ_d -J_d F_k^T)    + \braces{\frac{m}{h}(x_{k+1}-x_k) \hat d R_k^T e_1+\frac{m}{h}(y_{k+1}-y_k) \hat d R_k^T e_2%
    - \frac{h}{2} m g \hat d R_k^T e_3}^\wedge,\label{eqn:pwk}\\
R_{k+1} & =R_kF_k,\label{eqn:Rkp3}\\
p_{x_{k+1}} & = p_{x_k}+hu_{x_{k+1}},\label{eqn:pxkp}\\
p_{y_{k+1}} & = p_{y_k}+hu_{y_{k+1}},\label{eqn:pykp}\\
\hat p_{\Omega_{k+1}} & = \frac{1}{h} (J_dF_k -F_k^T J_d )    + \braces{\frac{m}{h}(x_{k+1}-x_k) \hat d R_{k+1}^T e_1%
    +\frac{m}{h}(y_{k+1}-y_k) \hat d R_{k+1}^T e_2%
    + \frac{h}{2} m g \hat d R_{k+1}^T e_3}^\wedge.\label{eqn:pwkp}
\end{align}}
The momenta variables $p_\Omega\in\Re^3$, $p_x,p_y\in\Re$ are given by
\begin{align}
    \begin{bmatrix}
    p_\Omega \\ p_x \\ p_y
    \end{bmatrix}
=
    \begin{bmatrix}
    J & m\hat d R^T e_1  & m\hat d R^T e_2\\
    -me_1^T R \hat d &  M+m & 0\\
    -me_2^T R \hat d & 0 & M+m
    \end{bmatrix}
    \begin{bmatrix}
    \Omega \\ \dot x \\ \dot y
    \end{bmatrix}
.\label{eqn:p}
\end{align}
The detailed derivation of this Lie group variational integrator is available in~\cite{Lee.2008}. For given $(R_k,x_k,y_k,\Omega_k,\dot x_k,\dot y_k)$, we compute $(p_{\Omega_k},p_{x_k},p_{y_k})$ by \refeqn{p}. We use a fixed-point iteration to compute $R_{k+1}$. For an initial guess for $R_{k+1}$, the corresponding $x_{k+1},y_{k+1}$ are obtained by using \refeqn{pxk},\refeqn{pyk}. Then, we can find $F_k$ by solving \refeqn{pwk}. The updated value for $R_{k+1}$ is given by \refeqn{Rkp3}. This is repeated until $R_{k+1}$ converges. Then, $x_{k+1},y_{k+1}$ are obtained from \refeqn{pxk},\refeqn{pyk}, and $(p_{\Omega_{k+1}},p_{x_{k+1}},p_{y_{k+1}})$ are obtained by \refeqn{pxkp},\refeqn{pykp}, and \refeqn{pwkp}. The velocities $(\Omega_{k+1},\dot x_{k+1},\dot y_{k+1})$ are obtained from \refeqn{p}. This yields a flow map,
\[(R_k,x_k,y_k,\Omega_k,\dot x_k,\dot y_k)\mapsto(R_{k+1},x_{k+1},y_{k+1},\Omega_{k+1},\dot x_{k+1},\dot y_{k+1}).\]

\subsubsection*{Optimal control problem}
The objective of the optimal control problem is to transfer the 3D pendulum on a cart from a given initial condition $(R_0,x_0,y_0,\Omega_0,\dot x_0,\dot y_0)$ to a desired terminal condition $(R^f,x^f,y^f,\Omega^f,\dot x^f,\dot y^f)$ during a fixed maneuver time $Nh$, while minimizing the square of the $l_2$ norm of the control inputs.
\begin{gather}
\min_{u_{k+1}} \braces{\mathcal{J}=\sum_{k=0}^{N-1}
\frac{h}{2}u_{k+1}^T W u_{k+1}},\label{eqn:J3dpendcart}
\end{gather}
where $u_{k}=[u_{x_k};u_{y_k}]\in\Re^2$, and $W\in\Re^{2\times 2}$ is a symmetric positive-definite matrix. The 3D pendulum on a cart is underactuated, since only the planar motion of the cart in its horizontal plane is actuated.

\subsubsection*{Computational approach}

We apply a direct optimal control approach. The control inputs are parameterized by several points that are uniformly distributed over the maneuver time, and control inputs between these points are approximated using a cubic spline interpolation. For given control input parameters, the value of the cost is given by \refeqn{J3dpendcart}, and the terminal conditions are obtained by the discrete-time equations of motion given by \refeqn{pxk}-\refeqn{pwkp}. The control input parameters are optimized using constrained nonlinear parameter optimization to satisfy the terminal boundary conditions while minimizing the cost.

This approach is computationally efficient when compared to the usual collocation methods, where the continuous-time equations of motion are imposed as constraints at a set of collocation points. Using the proposed discrete-time optimal control approach, optimal control inputs can be obtained by using a large step size, thereby resulting in efficient total computations. Since the computed optimal trajectories do not have numerical dissipation caused by conventional numerical integration schemes, they are numerically more robust. Furthermore, the corresponding gradient information is accurately computed, which improves the convergence properties of the numerical optimization procedure.

\subsubsection*{Numerical example}

Properties of the 3D pendulum and the cart are chosen as,
\[ M=m=1\,\mathrm{kg},\quad J=\mathrm{diag}[1.03,1.04,0.03]\,\mathrm{kgm^2},\quad\text{and}\quad d=[0,0,1]\,\mathrm{m}.\]
The desired maneuver is a rest-to-rest $180^\circ$ rotation of the pendulum about the vertical axis, while the cart returns to the initial location at the terminal time. The corresponding boundary conditions are given by
\begin{gather*}
    R_0=I,\quad \Omega_0=[0,0,0],\quad x_0=y_0=0,\quad \dot x_0=\dot y_0=0,\\
    R^f=\mathrm{diag}[-1,-1,1],\quad \Omega^f=[0,0,0],\quad x^f=y^f=0,\quad \dot x^f=\dot y^f=0.
\end{gather*}
The maneuver time is $2$ seconds, and the time step is $h=0.01$. Since only the planar motion of the cart is actuated, the rotation of the 3D pendulum is caused by the nonlinear coupling between the cart and the pendulum.

\begin{figure}[h]
\centerline{
    \subfigure[Optimal maneuver of a 3D pendulum on a cart]{
        \includegraphics[width=0.95\textwidth]{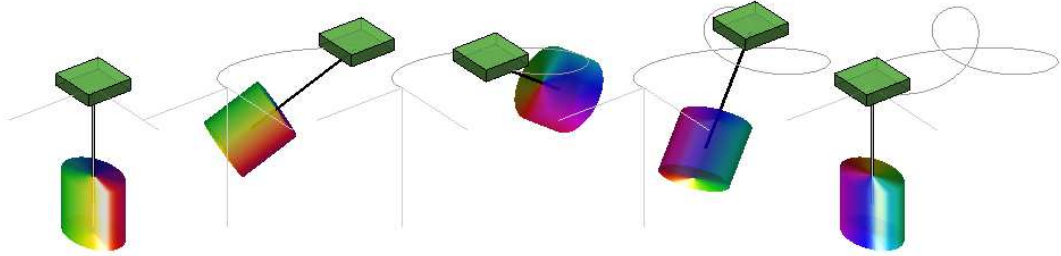}\label{fig:Opt3DPC_3d}}
}
\centerline{
    \subfigure[Control force $u=(u_x,u_y)$]{
        \includegraphics[width=0.46\textwidth]{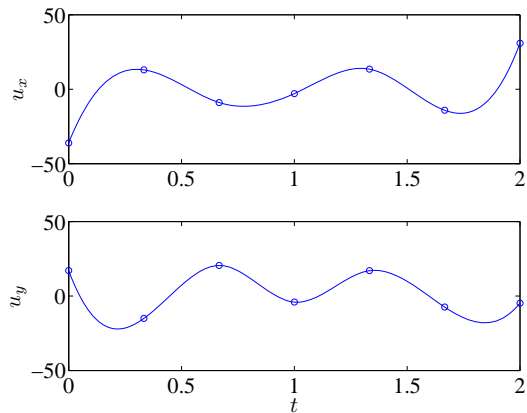}}
    \hfill
    \subfigure[Angular velocity $\Omega$]{
        \includegraphics[width=0.45\textwidth]{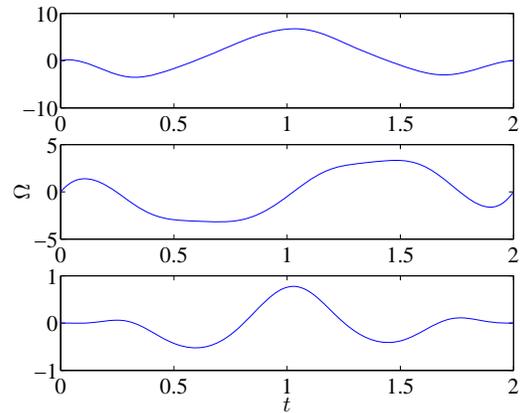}}
    }
\caption{Optimal control of a 3D pendulum on a cart}\label{fig:Opt3DPC}
\end{figure}

Each component of the control inputs is parameterized by 7 points. The resulting 14 control input parameters are optimized using sequential quadratic programming. \reffig{Opt3DPC} illustrates the optimal maneuver of the pendulum and the cart, angular velocity, and optimal control inputs. The blue circles denote the optimized control input parameters. The optimal cost and the violation of the terminal boundary conditions are
$297.43$, and $1.83\times10^{-8}$, respectively. The optimal motion of the cart on the horizontal plane consists of a triangular-shaped loop, and the optimal maneuver of the 3D pendulum consists of large angle rotations. This also demonstrates the advantages of the computational geometric optimal control approach: it is difficult to study this kind of aggressive maneuvers of a multibody system using local coordinates, due to the coordinate singularities and the complexity of the equations in local coordinates. The presented computational geometric optimal control approach accurately characterizes the nonlinear coupling between the cart and the pendulum dynamics to obtain a nontrivial optimal maneuver of the 3D pendulum on a cart.

\subsection{Optimal attitude reorientation of two connected rigid bodies on $\SO\times\SO$}

Consider two rigid bodies connected with a ball joint that has three rotational degrees of freedom. This represents a freely rotating system of coupled rigid bodies. The relative equilibria structure of this rigid body dynamics has been studied in \cite{Wan.1990}. We introduce three frames; an inertial frame and two body-fixed frames. Define
\begin{center}
\begin{tabular}{lp{13cm}}
$x\in\Re^3$ & Position of the ball joint in a reference frame\\
$R_i\in\SO$ & Rotation matrix from the $i$-th body-fixed frame to a reference frame\\
$d_i\in\Re^3$ & Vector from the joint to the mass center of the $i$-th body in the $i$-th body-fixed frame\\
$m_i\in\Re$ & Mass of the $i$-th body
\end{tabular}
\end{center}
for $i\in\{1,2\}$. The configuration manifold is $\SO\times\SO\times\Re^3$. In the absence of the potential field, the connected rigid body model has two symmetries; a symmetry of the translational action of $\Re^3$, and a symmetry of the rotational action of $\SO$.\footnote{These can be considered as a single symmetry of the translational and rotational action of $\SE$, but they are considered separately in this optimal control problem. By the general theory of reduction by stages~\cite{CeMaRa2001}, the two approaches are equivalent.} Due to these symmetries, the total linear momentum and the total angular momentum are preserved, and the configuration manifold can be reduced to a quotient space.

In this optimal control problem, we reduce the configuration manifold to $\SO\times\SO$ using the symmetry of the translational action of $\Re^3$. The corresponding value of the total linear momentum is set to zero. The resulting connected rigid bodies model with a fixed mass center is closely related to the falling cat problem~\cite{Eno.BK1993}. An appropriate cyclic change in the shape of the body yields a rotation in the orientation of the cat in accordance with the geometric phase effect~\cite{Mon.1991}. In contrast to other models of the falling cat, which typically introduce two one-dimensional rotational joints, with a shape space given by $S^1\times S^1$, we consider instead a single ball joint with a shape space given by $SO(3)$.

Similar to the falling cat problem, we assume that an internal control moment $u\in\Re^3$ is applied at the joint, so that it controls the relative attitude between two rigid bodies. More precisely, the control input $u$ represents the control moment applied to the first rigid body, represented in the reference frame. The equal and opposite control moment is applied to the second rigid body. Therefore the control moment changes the shape of the system. The total angular momentum is conserved for the controlled dynamics as the control input is an internal moment of the connected rigid bodies system. This optimal control problem is similar to the optimal control problem of the 3D pendulum discussed in Section~\ref{subsec:3dpend}, as the control input respects the symmetry, and the corresponding momentum is preserved in the controlled dynamics.

The discrete Lagrangian for the two connected rigid bodies is
\begin{align}
L_d (R_{1_k},F_{1_k}, R_{2_k},F_{2_k},x_k,x_{k+1}) =& \frac{m_1+m_2}{2h}(x_{k+1}-x_k)\cdot(x_{k+1}-x_k) + \frac{1}{h}\tr{(I_{3\times  3}-F_{1_k})J_{d_1}}\nonumber\\
&+\frac{1}{h}\tr{(I_{3\times  3}-F_{2_k})J_{d_2}}
+\frac{1}{h}\tr{m_1R_{1_k}(F_{1_k}-I_{3\times 3})d_1(x_{k+1}-x_k)^T}\nonumber\\
&+\frac{1}{h}\tr{m_2R_{2_k}(F_{2_k}-I_{3\times 3})d_2(x_{k+1}-x_k)^T}.
\end{align}
From \refeqn{DEP}, we obtain the Lie group variational integrator, viewed as discrete-time equations of motion on $\SO\times\SO\times\Re^3$. Since we are only interested in rotational maneuvers, we derive the following reduced equations of motion on $\SO\times\SO$ using the fact that the linear momentum is conserved.
{\allowdisplaybreaks
\begin{align}
\hat p_{1_k}  = & \frac{1}{h} \braces{F_{1_k} (J_{d_1}-\alpha m_1 d_1 d_1^T) -  (J_{d_1}-\alpha m_1 d_1 d_1^T ) F_{1_k}^T}\nonumber\\
&-\beta\frac{m_1}{h}( R_{1_k}^T R_{2_k} F_{2_k} d_2 d_1^T - d_1 d_2^T F_{2_k}^T R_{2_k}^T R_{1_k})
+ \beta\frac{m_1}{h} (R_{1_k}^T R_{2_k} d_2 d_1^T -d_1 d_2^T R_{2_k}^T R_{1_k}) ,\label{eqn:p1kd}\\
\hat p_{2_k}  = & \frac{1}{h}\braces{ F_{2_k} (J_{d_2} -\beta m_2 d_2d_2^T) -(J_{d_2}-\beta m_2 d_2d_2^T)F_{2_k}^T}\nonumber\\
&-\alpha\frac{m_2}{h} ( R_{2_k}^T R_{1_k} F_{1_k} d_1 d_2^T - d_2 d_1^T F_{1_k}^T R_{1_k}^T R_{2_k}) +\alpha\frac{m_2}{h} (R_{2_k}^T R_{1_k} d_1 d_2^T -d_2 d_1^T R_{1_k}^T R_{2_k}) ,\label{eqn:p2kd}\\
R_{1_{k+1}}= &R_{1_k}F_{1_k},\label{eqn:R1kp}\\
R_{2_{k+1}}= &R_{2_k}F_{2_k},\label{eqn:R2kp}\\
p_{1_{k+1}} = & F_{1_k}^T (p_{1_k} -(B_{1_k}-B_{1_k}^T)^\vee)+hR_{1_{k+1}}^T u_{k+1},\label{eqn:p1kp}\\
p_{2_{k+1}} = & F_{2_k}^T (p_{2_k} -(B_{2_k}-B_{2_k}^T)^\vee)-hR_{2_{k+1}}^T u_{k+1},\label{eqn:p2kp}
\end{align}}
where $\alpha=\frac{m_1}{m_1+m_2}$, $\beta=\frac{m_2}{m_1+m_2}\in\Re$, and the matrix $B_{i_k}\in\Re^{3\times 3}$ for $i\in\{1,2\}$ is defined as
\begin{align}
B_{i_k} = \frac{m_i}{h}(F_{i_k}-I)d_i\braces{-\alpha R_{1_k} (F_{1_k}-I)d_1 -\beta R_{2_k} (F_{2_k}-I)d_2}^T R_{i_k}.\label{eqn:Bik}
\end{align}
The momenta variables $p_{1},p_{2}\in\Re^3$ are given by
\begin{align}
\begin{bmatrix} p_{1} \\
    p_{2} \end{bmatrix}=
    \begin{bmatrix} J_1+\alpha m_1\hat d_1^T & \beta m_1 \hat d_1 R_{1}^T R_{2}\hat d_2\\
    \alpha m_2\hat d_2 R_{2}^T R_{1} \hat d_1 & J_2 + \beta m_2 \hat d_2^2\end{bmatrix}
    \begin{bmatrix} \Omega_{1} \\ \Omega_{2} \end{bmatrix}.\label{eqn:LTp}
\end{align}

For given $(R_{1_k},R_{2_k},\Omega_{1_k},\Omega_{2_k})$, we find $p_{1_k},p_{2_k}$ by \refeqn{LTp}. We solve the implicit equations \refeqn{p1kd}, \refeqn{p2kd} to obtain $F_{1_k},F_{2_k}$. Then, $R_{1_{k+1}},R_{2_{k+1}}$ are obtained from \refeqn{R1kp},\refeqn{R2kp}, and $p_{1_{k+1}},p_{2_{k+1}}$ are obtained by \refeqn{p1kp},\refeqn{p2kp}. Finally, $\Omega_{1_{k+1}},\Omega_{2_{k+1}}$ are computed from \refeqn{LTp}. This yields a discrete flow map $(R_{1_k},R_{2_k},\Omega_{1_k},\Omega_{2_k})\mapsto(R_{1_{k+1}},R_{2_{k+1}},\Omega_{1_{k+1}},\Omega_{2_{k+1}})$.

\subsubsection*{Optimal control problem}
The objective of the optimal control problem is to transfer the connected rigid bodies from a given initial condition $(R_{1_0},R_{2_0},\Omega_{1_0},\Omega_{2_0})$ to a desired terminal condition $(R_1^f,R_2^f,\Omega_1^f,\Omega_2^f)$ during a fixed maneuver time $Nh$, while minimizing the square of the $l_2$ norm of the control inputs.
\begin{gather}
\min_{u_{k+1}} \braces{\mathcal{J}=\sum_{k=0}^{N-1}
\frac{h}{2}u_{k+1}^T W u_{k+1}},\label{eqn:J4}
\end{gather}
where $W\in\Re^{3\times 3}$ is a symmetric positive-definite matrix. In particular, we choose an attitude maneuver that is described by a rest-to-rest rotation of the entire system while the relative attitude configuration at the terminal time is the same as at the initial time.

\subsubsection*{Computational approach}

We apply a direct optimal control approach. For a given control input, the value of the cost is given by \refeqn{J4}, and the terminal conditions are obtained by the discrete-time equations of motion given by \refeqn{p1kd}-\refeqn{Bik}. We use constrained nonlinear parameter optimization to minimize the cost function subject to the terminal boundary condition obtained by the discrete-time equations of motion.

Since the total angular momentum is conserved regardless of the control input, the terminal constraints introduces a singularity due to the rotational symmetry. This ill-conditioning can be avoided by disregarding the terminal angular velocity constraint for the second body. For the given boundary conditions, the terminal angular velocity condition is automatically satisfied if the remaining terminal constraints are satisfied, due to the angular momentum conservation property. By formulating the optimization process this way, we eliminate the source of numerical ill-conditioning. This is similar to the modified computational approach discussed in Section~\ref{subsec:3dpend}.

\subsubsection*{Numerical example}

Properties of the rigid bodies are chosen as
\begin{align*}
    m_1 = 1.5\mathrm{kg},\quad J_1=\begin{bmatrix}0.18 & 0.32 & 0.32\\0.32 & 1.88 & -0.06\\0.32 & -0.06& 1.86\end{bmatrix}\mathrm{kg\cdot m^2},\quad d_1=[-1.08,0.20,0.20]\mathrm{m},\\
    m_2 = 1\mathrm{kg},\quad J_2=\begin{bmatrix}0.11 &-0.18 &-0.18\\-0.18 & 0.89 & -0.04\\-0.18 & -0.04 & 0.88\end{bmatrix}\mathrm{kg\cdot m^2},\quad d_2=[0.9,0.2,0.2]\mathrm{m}.
\end{align*}
The desired maneuver is a rest-to-rest $180^\circ$ rotation about the $x$ axis.
\begin{gather*}
    R_{1_0}=I,\quad \Omega_{1_0}=0,\quad     R_{2_0}=I,\quad \Omega_{2_0}=0,\\
    R_1^f=\mathrm{diag}[1,-1,-1] ,\quad \Omega_1^f=[0,0,0],\quad     R_2^f=\mathrm{diag}[1,-1,-1],\quad \Omega_2^f=[0,0,0].
\end{gather*}
The maneuver time is $4$ seconds, and the step size is $h=0.01$.

\begin{figure}[htb]
\centerline{
    \subfigure[Optimal maneuver]{
        \includegraphics[width=0.85\textwidth]{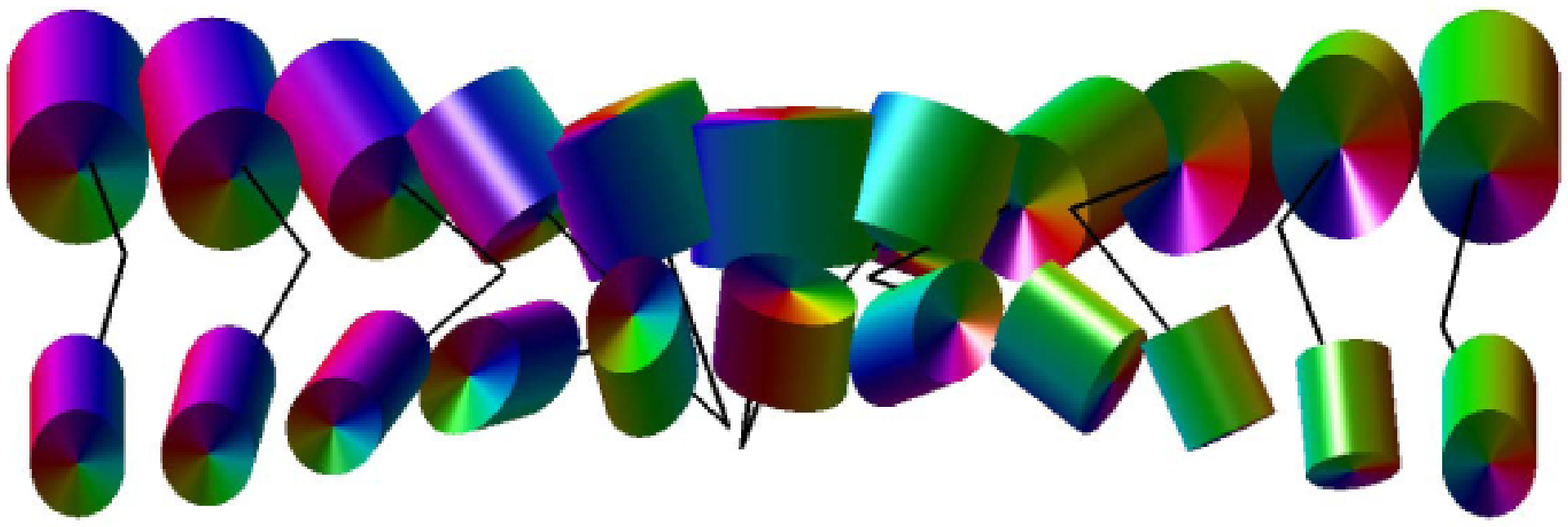}\label{fig:OptTRB_3d}}
}
\centerline{
    \subfigure[Control input $u$]{
        \includegraphics[width=0.45\textwidth]{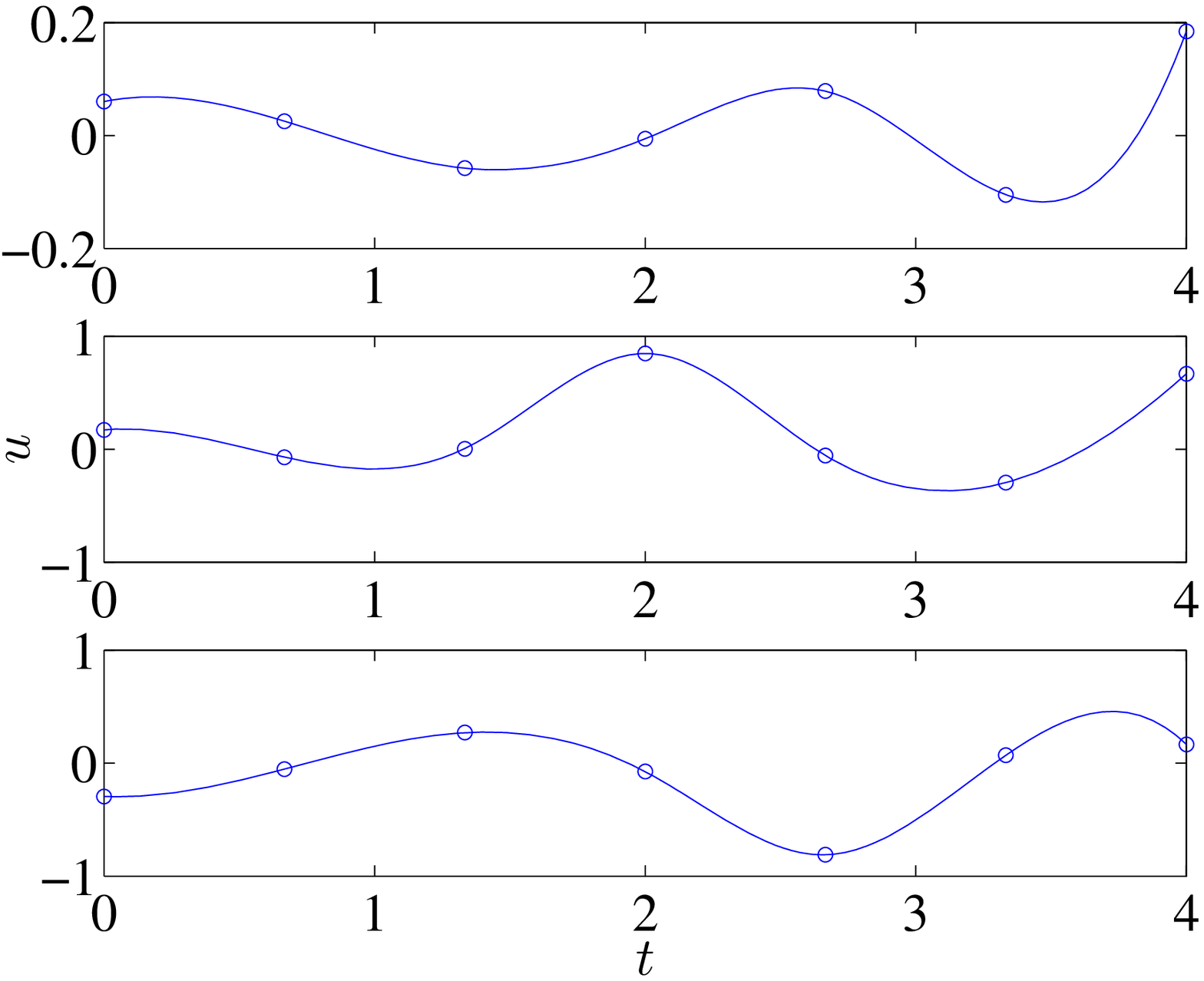}}
    \hfill
    \subfigure[Angular velocity ($\Omega_1$:solid, $\Omega_2$:dashed)]{
        \includegraphics[width=0.45\textwidth]{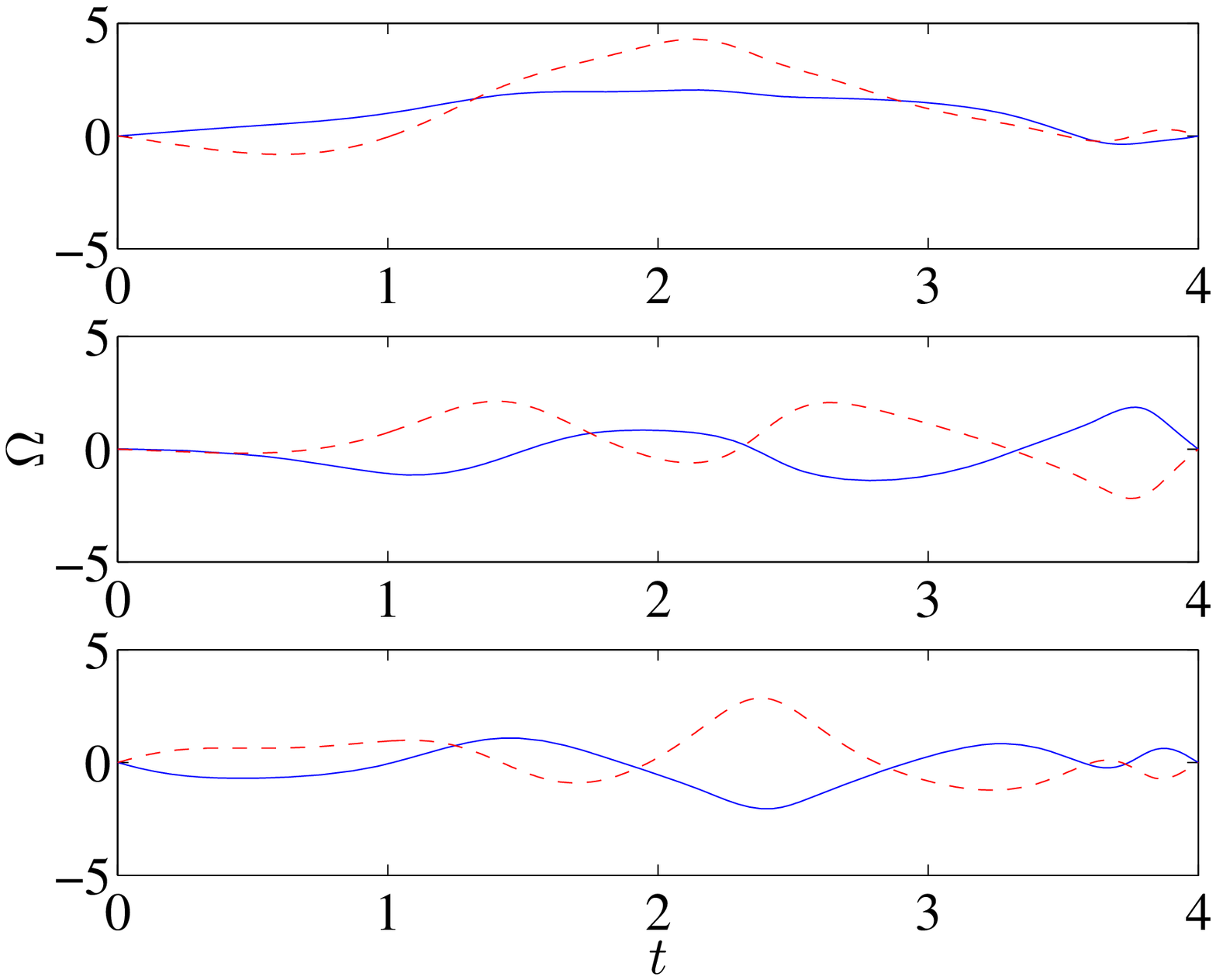}}
}
\caption{Optimal control of two connected rigid bodies}\label{fig:OptTRB}
\end{figure}

We parameterize each component of the control input at $7$ discrete points, and the control inputs are reconstructed by cubic spline interpolation. The resulting $21$ control input parameters are optimized by a sequential quadratic programming method to satisfy the terminal boundary conditions while minimizing the cost function.

\reffig{OptTRB} shows the optimal maneuver of the rigid bodies, angular velocity, and optimal control inputs. The blue circles denote the optimized control input parameters. The optimal cost and the violation of the terminal boundary conditions are
$0.574$, and $2.48\times10^{-8}$, respectively.

The optimal maneuver consists of large angle rotations of the two rigid bodies. Throughout this complicated maneuver, the total angular momentum is zero, and the rotation about the $e_1$ axis  depends on the geometric phase effect. This also demonstrates the advantages of the computational geometric optimal control approach. The Lie group variational integrator computes the weak geometric phase effect accurately, so that the iterations converge to a nontrivial optimal maneuver of the coupled rigid bodies.

\section{Conclusions}

In this paper, a computational geometric approach for an optimal control problem of rigid body dynamics has been developed. The essential idea is formulating a discrete-time optimal control problem using a structure-preserving geometric numerical integrator, referred to as a Lie group variational integrator, and applying standard optimal control approaches, such as an indirect optimal control and a direct optimal control, to discrete-time equations of motion. This method is in contrast to the usual optimal control approach, where the discretization appears only in the last stage when numerically computing the optimal control inputs.

The computational geometric optimal control approach has substantial advantages in terms of preserving the geometric properties of optimality conditions. The discrete flow of Lie group variational integrators has desirable geometric properties, such as symplecticity and momentum preservation, and it is more reliable and robust over longer time periods. The computational geometric optimal control approach inherits the desirable properties of the Lie group variational integrator. In the necessary conditions for optimality, the multiplier equations are dual to the linearized equations of motion. Since the linearized flow of a Lagrangian/Hamiltonian system is also symplectic, the multiplier equations inherit certain geometric properties. The discrete-time necessary conditions preserve the geometric properties of the optimality conditions, as they are derived from a discrete-time analogue of Hamilton's variational principle that yields a symplectic discrete-time flow.

The computational geometric optimal control approach allows us to find the optimal control input more efficiently. In the indirect optimal control, the shooting method may be prone to numerical ill-conditioning, as a small change in the initial multiplier can cause highly nonlinear behavior of the terminal condition. However, as shown in \reffig{OptRB_err} and \reffig{Opt3DP_err}, the computational geometric optimal control approach exhibits excellent numerical convergence properties. This is because the proposed computational algorithms are geometrically exact and numerically accurate. There is no numerical dissipation introduced by the numerical algorithm, and therefore, we are more accurately characterizing the sensitivities along the solution.

Another advantage of the computational geometric optimal control of rigid bodies is that the method is directly developed on a Lie group. There is no ambiguity or singularity in representing the configuration of rigid bodies globally, and the resulting equations of motion are more compact than those written in terms of local coordinates. As illustrated by \reffig{Opt3DPC_3d} and \reffig{OptTRB_3d}, the presented computational geometric optimal control approach utilizes the effects of the nonlinear coupling and the weak geometric phase of a multibody system to obtain nontrivial aggressive maneuvers of the rigid bodies. These results are independent of a specific choice of local coordinates, and they completely avoid any singularity, ambiguity, and complexity associated with local coordinates. Furthermore, the numerical results are group-equivariant, and are independent of the choice of inertial frame, which is in contrast to methods based on local coordinate representations. By formulating the problem in a global and intrinsic fashion, the algorithms presented are able to explore the space of control strategies which extend beyond a single coordinate chart, thereby providing a deeper insight into the global controlled dynamics of systems of rigid bodies.

\bibliography{cis,tylee}
\bibliographystyle{IEEEtranS}

\end{document}